\documentclass[a4paper,12pt]{article}
\usepackage{amsmath,amssymb,amsthm,latexsym,graphicx}
\usepackage{marvosym}
\usepackage{dictsym}
\usepackage{natbib}
\usepackage{longtable}
\usepackage{hyperref}
\usepackage{makeidx}

\usepackage{a4wide}

\hypersetup{
	colorlinks=true,
    breaklinks=true,
    filecolor=blue,
    urlcolor=blue,
    linkcolor=blue,
    citecolor=blue,
    anchorcolor=blue,
    frenchlinks=true,
    pdfborder={0 0 0},
    pdfdisplaydoctitle=true,
    bookmarksnumbered=true
}

\theoremstyle{definition}

\newtheorem{lemma}{Lemma}
\newtheorem{proposition}[lemma]{Proposition}
\newtheorem{theorem}[lemma]{Theorem}

\newtheorem{definition}[lemma]{Definition}
\newtheorem{example}[lemma]{Example}

\newtheorem{tab}[lemma]{Table}
\newtheorem{algorithm}[lemma]{Algorithm}

\renewenvironment{proof}{\textbf{Proof:}}{\hfill\Coffeecup}

\makeindex
\title{Families of bitangent planes of space curves and minimal non-fibration families}
\author{Niels Lubbes}

\date{\today}

\begin{document}\maketitle
\begin{abstract}
 We define a cone curve to be a reduced sextic space curve which lies on a quadric cone and does not go through the vertex. We classify families of bitangent planes of cone curves. The methods we apply can be used for any space curve with ADE singularities, though in this paper we concentrate on cone curves.

 An embedded complex projective surface which is adjoint to a degree one weak Del Pezzo surface contains families of minimal degree rational curves, which cannot be defined by the fibers of a map. Such families are called minimal non-fibration families. Families of bitangent planes of cone curves correspond to minimal non-fibration families. The main motivation of this paper is to classify minimal non-fibration families.

 We present algorithms wich compute all bitangent families of a given cone curve and their geometric genera. We consider cone curves to be equivalent if they have the same singularity configuration. For each equivalence class of cone curves we determine the possible number of bitangent families and the number of rational bitangent families. Finally we compute an example of a minimal non-fibration family on an embedded weak degree one Del Pezzo surface.
\end{abstract}

\newpage
\tableofcontents
\newpage

\section{Introduction}

\subsection{Problems}
A cone curve ${\textrm{C}}$ is a reduced sextic curve which lies on the quadric cone ${\textrm{Q}}$ in three-space. We also assume that a cone curve does not go through the vertex of ${\textrm{Q}}$. Cone curves occur as branching curves of 2:1 coverings of ${\textrm{Q}}$ by weak degree one Del Pezzo surfaces. A family of bitangent planes $F$ of ${\textrm{C}}$ (bitangent family for short) can be defined as an irreducible component of $U$ where \[ U=\{~(a,b)~~|~~ \textrm{there exists a plane which is bitangent at } a {\textrm{~and~}} b ~\} \subset {\textrm{C}}\times {\textrm{C}}. \] We call a bitangent family $F$ rational if the geometric genus $p_g F$ is zero.

 In this paper we concentrate on cone curves, but many of our methods can be used to find bitangent families of arbitrary space curves.

 In this paper we present a solution to the following problem:

 \textbf{Problem 1.} Compute all bitangent families of a given cone curve and their geometric genera.

 The solution to this problem is presented by Algorithm~\ref{alg:f2_dual_tangent_developable} and Algorithm~\ref{alg:f2_analyze_developable_surface}. However, the algorithm may not terminate within reasonable time because of elimination algorithms. Algorithm 120 in \cite{nls2}, chapter 8, section 11, page 172, does not use elimination and computes the geometric genera of the bitangent families $F$. However, again the algorithm might not terminate, because of the complexity of polynomial expansion.

 We will define two cone curves to be equivalent if they have the same ADE singularity configuration. Each equivalence class can be represented by a root subsystem of the root system with Dynkin type $E_8$. We shall represent such root subsystems by so called C1 labels. We recall the definitions in the next section and the details can be found in \cite{nls-f1}. Up to equivalence there is a finite list of cone curves. We will see that the degree of the components of a cone curve is uniquely defined for each equivalence class. For each entry we consider the bitangent families and their geometric genera. In this paper we present a solution to the following problem:

 \textbf{Problem 2.} For each equivalence class of cone curves determine the possible number of bitangent families and the number of rational bitangent families.

 The solution of this problem is presented at Theorem~\ref{thm:f2_coc_cls}. It turns out that for some equivalence classes of cone curves not all representatives have the same number of bitangent families (see for example Proposition~\ref{prop:f2_cone_curve_genus_two}). In this case an upper and a lower bound is presented. It is remarkable that a cone curve with $4A_2$ singularities is the unique equivalence class which has no family of bitangent planes, other then the trivial family of tritangent planes which is defined by the ruling of the quadric cone ${\textrm{Q}}$.

 We consider the bitangent families of generic cone curves in more detail:

 \textbf{Problem 3.} Determine the arithmetic genera of the bitangent families of a generic cone curve. Moreover, find bounds on the number of special planes in the bitangent families.

 The solution to this problem is presented by Theorem~\ref{thm:f2_cone_curve_cor}.

 We shall see that ${\textrm{C}}$ defines up to projective isomorphism the anticanonical model of a weak degree one Del Pezzo surface ${\textrm{Z}}$. A bitangent family $F$ which is not defined by the ruling of ${\textrm{Q}}$ defines a minimal non-fibration family of ${\textrm{Z}}$. We call such families T5-families (see chapter 7, section 1, definition 85 in \cite{nls2}).

 \textbf{Problem 4.} Compute the T5 families of a given weak degree one Del Pezzo surface.

 We propose a solution towards this problem in terms of an example: Example~\ref{ex:f2_T5}.

\subsection{Motivation}
We recall that a weak Del Pezzo surface has a nef and big anticanonical divisor class $-K$. Let ${\textrm{X}}$ be a weak Del Pezzo surface of degree $K^2=1$. The anticanonical model of ${\textrm{X}}$ is $\varphi_{-3K}({\textrm{X}})\subset{\textbf{P}}^6$ where $\varphi_{-3K}$ is the map associated to the class $-3K$. We have that $\ensuremath{{\textrm{X}}\stackrel{\varphi_{-2K}}{\rightarrow}{\textrm{Q}}}$ defines a 2:1 covering of the quadric cone ${\textrm{Q}}$. We can pull back a bitangent family along this covering. If the bitangent family is not the tritangent planes defined by the ruling of ${\textrm{Q}}$ then this pull back defines a family of minimal degree rational curves on ${\textrm{X}} \subset {\textbf{P}}^6$. This family of curves cannot be defined by the fibres of a morphism. We call such families minimal non-fibration families.

 The minimal non-fibration families which are rational ($p_g F=0$) correspond to unirational parametrizations. In order to illustrate this let us consider an example of a minimal non-fibration family in the complex plane. Although the plane is a degree nine instead of degree one Del Pezzo surface, the following example illustrates the idea. For the analogues example for a weak degree one Del Pezzo surface see Example~\ref{ex:f2_T5}. A minimal family $U \subset C \times {\textbf{C}}^2$ of lines in the complex plane tangent to the unit circle $C$ is defined as follows: $C: a^2+b^2-1=0$ and $U: ax+by-1=0$.
\begin{center}
 {\includegraphics[width=3cm,height=3cm]{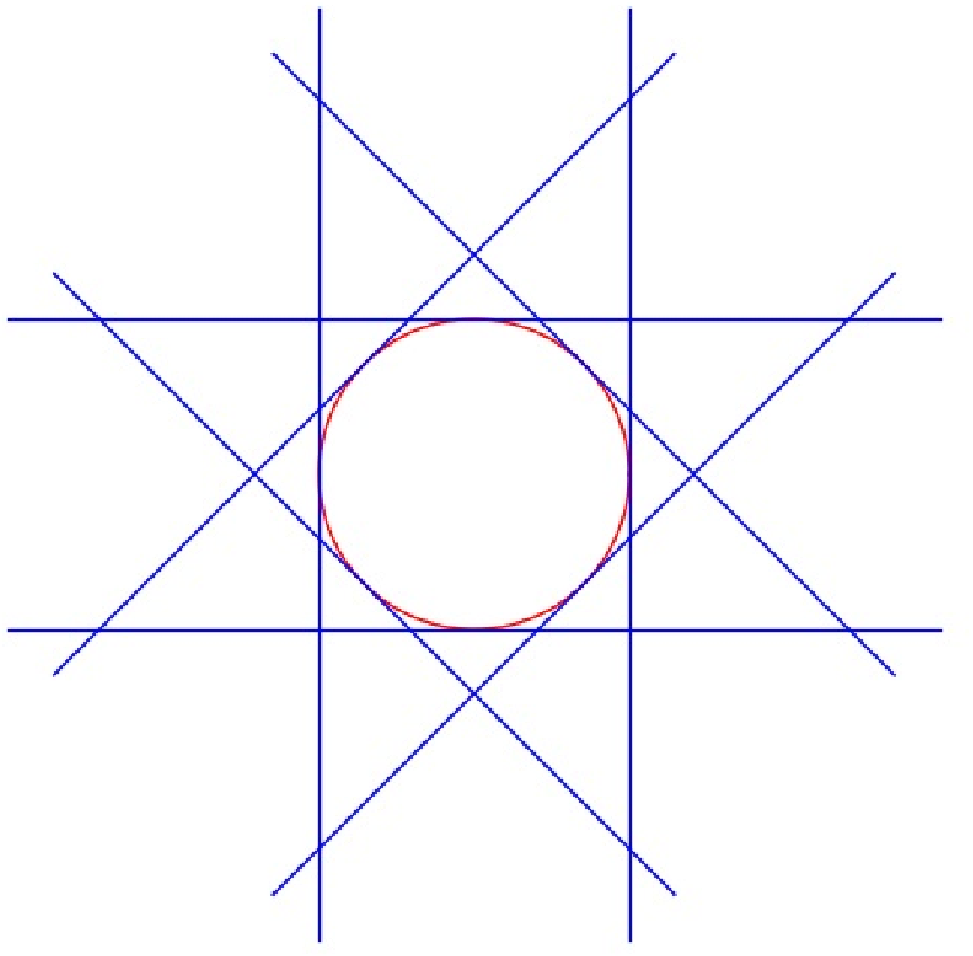}}
\end{center} We have that $s \mapsto (~f(s),~g(s)~):=(~\frac{1-s^2}{1+s^2},~ \frac{2s}{1+s^2} ~)$ is a parametrization of $C$. It follows that $(s,t)\mapsto (~f(s),~g(s)~;~t,~\frac{1-f(s)t}{g(s)}~)$ parametrizes $U$. From the figure above we see that every point in the plane is reached by two lines tangent to the unit circle. It follows that the second projection map $\ensuremath{U\stackrel{}{\rightarrow}{\textbf{C}}^2}$ is a two to one map. We find that the composition of the parametrization of $U$ with the projection $(s,t) \mapsto (~t,~\frac{1-f(s)t}{g(s)}~)$ is a unirational parametrization which is of minimal degree with respect to $s$. If we fix $s$ we parametrize a line in the plane. Note that the degree  with respect to  $t$ depends on the degree of $C$.

 Minimal families are defined as families of minimal degree rational curves. This paper is part of the larger classification of minimal families on surfaces. In \cite{nls1} the classification of minimal families on surfaces is reduced to the classification of minimal families on geometrically ruled surfaces and weak Del Pezzo surfaces. Of these two types of surfaces only weak Del Pezzo surfaces of degree one, two or nine have non-fibration families which are minimal.

 The classification of minimal non-fibration families of degree two Del Pezzo surfaces are determined in an analogues but less involved way as for degree one Del Pezzo surfaces. A weak degree two Del Pezzo surface ${\textrm{S}}$ admits a 2:1 cover of the projective plane, with a quartic plane curve ${\textrm{B}}$ as branching curve (see for example \cite{stu2}). The family of tangent lines of the quartic plane curve are determined by the non-linear components of ${\textrm{B}}$. The families of tangent lines pull back along the 2:1 covering to minimal non-fibration families on ${\textrm{S}}$. The details can be found in \cite{nls2}.

 So one motivation of this paper can be stated as follows: Classification of minimal non-fibration families on surfaces.

 The algebraic methods in this paper apply to arbitrary space curves. We use techniques from \cite{jha2}, chapter 2, section 5, to state properties about bitangent families. In this paper we consider all ADE singularities, and not just the traditional ones, which we hope to be instructive.

 In \cite{stu1} multi-tangent planes of space curves are considered as an application to compute the convex hull of space curves.

 \newpage
\subsection{Overview}
We start by recalling some theory concerning degree one weak Del Pezzo surfaces. We use this theory since cone curves occur as branching curves of 2:1 coverings of the quadric cone by degree one weak Del Pezzo surfaces.

 In the next section we recall the `C1 classification' which is the classification of the isomorphism classes of root subsystems of $E_8$. Each root subsystem will be represented by a C1 label. In our setting a C1 label represents a set of effective minus two classes in the Picard group of a weak degree one Del Pezzo surface.

 Next we consider some properties of cone curves: its canonical divisor class, genera and linear series of degree two and dimension one. We also formalize the relation between cone curves and weak degree one Del Pezzo surfaces.

 The C1 classification classifies the singularities of cone curves. For each C1 label in the C1 classification we determine the components of a cone curve with corresponding singularity configuration. Moreover we determine the geometric genus of each of the components.

 In the germs section we recall some definitions and the classification of space germs. We need this for the next section, where we consider the ramifications of projections of curve curves.

 After that we analyze the reduced bitangent correspondence on generic cone curves. We determine the arithmetic genus of the correspondence, and bound the number of special hyperplanes in the reduced bitangent correspondence.

 We consider the tangent developables of dual cone curves and the tangent developable of cone curves. The singular locus of these tangent developables determine bitangent families.

 In the following sections we consider bitangent families of  respectively  irreducible and reducible cone curves. All the results obtained are put together in Theorem~\ref{thm:f2_coc_cls} which states a table with the number of (rational) bitangent families of cone curves for each singularity configuration.

 The remaining sections discuss algorithms for determining the bitangent families of a given cone curve. We consider an elliptic cone curve with three cusps as a running example. We give an example of a minimal non-fibration family of a degree one Del Pezzo surface.

\newpage
{\Large\textbf{Overview table of environments}}

\newcounter{Movct}
\setcounter{Movct}{0}
\newcommand{\MovctAdd}{\addtocounter{Movct}{1}\arabic{Movct}}


\newcounter{Msection}
\setcounter{Msection}{0}
\newcommand{\MsectionAdd}{\addtocounter{Msection}{1}\setcounter{Msubsection}{0}\arabic{Msection}}

\newcounter{Msubsection}
\setcounter{Msubsection}{0}
\newcommand{\MsubsectionAdd}{\addtocounter{Msubsection}{1}\arabic{Msection}.\arabic{Msubsection}}
{\small
\begin{longtable}{|l|l|l|}
\hline
type  & \# & description \\\hline
\hline
\hline
\textbf{Section} &\MsectionAdd  & Introduction \\\hline\hline
\hline
\textbf{Subsection} &\MsubsectionAdd  & Problems \\\hline\hline
\hline
\textbf{Subsection} &\MsubsectionAdd  & Motivation \\\hline\hline
\hline
\textbf{Subsection} &\MsubsectionAdd  & Overview \\\hline\hline
\hline
\textbf{Subsection} &\MsubsectionAdd  & Notation \\\hline\hline
\hline
\textbf{Section} &\MsectionAdd  & Weak Del Pezzo surfaces of degree one \\\hline\hline

Definition &\MovctAdd &DP1 ring    \\\hline

Proposition &\MovctAdd &properties of DP1 ring    \\\hline

Proposition &\MovctAdd &canonical projection of degree one weak Del Pezzo surface    \\\hline
\hline
\textbf{Section} &\MsectionAdd  & C1 classification \\\hline\hline
\hline
\textbf{Section} &\MsectionAdd  & Properties of cone curves \\\hline\hline

Definition &\MovctAdd &cone curve    \\\hline

Proposition &\MovctAdd &properties of divisor classes of cone curve    \\\hline

Definition &\MovctAdd &weighted cone curve    \\\hline

Proposition &\MovctAdd &properties of weighted cone curve    \\\hline

Definition &\MovctAdd &cone curve to DP1 function    \\\hline

Proposition &\MovctAdd &properties of cone curve to DP1 function    \\\hline
\hline
\textbf{Section} &\MsectionAdd  & Components and singularities of cone curves \\\hline\hline

Proposition &\MovctAdd &properties of singularities of cone curve    \\\hline

Proposition &\MovctAdd &properties of components of cone curve    \\\hline
\hline
\textbf{Section} &\MsectionAdd  & Germs \\\hline\hline

Definition &\MovctAdd &ADE-singularities    \\\hline

Proposition &\MovctAdd &classification of simple space germs    \\\hline

Proposition &\MovctAdd &properties of space germs of space curves    \\\hline

Proposition &\MovctAdd &properties of delta invariant of ADE singularities    \\\hline
\hline
\textbf{Section} &\MsectionAdd  & Ramifications of projections of space curves \\\hline\hline

Proposition &\MovctAdd &ramification of projections of singularities    \\\hline
\hline
\textbf{Section} &\MsectionAdd  & Correspondences of cone curves \\\hline\hline

Definition &\MovctAdd &correspondences of curves    \\\hline

Definition &\MovctAdd &attributes of correspondences of curves    \\\hline

Proposition &\MovctAdd &properties of correspondences of curves    \\\hline

Definition &\MovctAdd &correspondences of cone curves    \\\hline

Theorem &\MovctAdd &properties of correspondences of cone curves    \\\hline
\hline
\textbf{Section} &\MsectionAdd  & Tangent developables of cone curves and its duals \\\hline\hline

Proposition &\MovctAdd &properties of tangent developable of cone curve and its dual    \\\hline
\hline
\textbf{Section} &\MsectionAdd  & Bitangent correspondences of irreducible cone curves \\\hline\hline

Proposition &\MovctAdd &properties of number of components of reduced bitangent correspondence: (6)    \\\hline

Proposition &\MovctAdd &properties of genus of components of reduced bitangent correspondence    \\\hline

Proposition &\MovctAdd &properties of bitangent correspondence of genus two cone curve    \\\hline
\hline
\textbf{Section} &\MsectionAdd  & Bitangent correspondences of reducible cone curves \\\hline\hline

Proposition &\MovctAdd &properties of bitangent correspondence: (2,2,2)    \\\hline

Proposition &\MovctAdd &properties of bitangent correspondence: (4,2)    \\\hline
\hline
\textbf{Section} &\MsectionAdd  & Classification of bitangent correspondence of cone curves \\\hline\hline

Theorem &\MovctAdd &properties of table of components of bitangent correspondence of cone curves    \\\hline

Table &\MovctAdd &table of components of bitangent correspondence of cone curves    \\\hline
\hline
\textbf{Section} &\MsectionAdd  & Algorithm: radical decomposition of equidimensional ideals \\\hline\hline

Algorithm &\MovctAdd &radical decomposition of equidimensional ideals    \\\hline

Proposition &\MovctAdd &radical decomposition of equidimensional ideals    \\\hline
\hline
\textbf{Section} &\MsectionAdd  & Algorithm: analyze tangent developable surface of dual of cone curve \\\hline\hline

Algorithm &\MovctAdd &tangent developable of dual of cone curve    \\\hline

Proposition &\MovctAdd &properties of analyze dual surface algorithm    \\\hline

Example &\MovctAdd &analyze dual surface algorithm    \\\hline
\hline
\textbf{Section} &\MsectionAdd  & Algorithm: analyze tangent developable surface of cone curve \\\hline\hline

Algorithm &\MovctAdd &analyze developable surface of cone curve    \\\hline

Proposition &\MovctAdd &properties of AnalyseDevelopableSurface algorithm    \\\hline

Example &\MovctAdd &analyze developable surface algorithm    \\\hline
\hline
\textbf{Section} &\MsectionAdd  & Example non-fibration family on a weak degree one Del Pezzo surface \\\hline\hline

Example &\MovctAdd &non-fibration family on a weak degree one Del Pezzo surface    \\\hline
\hline
\textbf{Section} &\MsectionAdd  & Acknowledgements \\\hline\hline
\end{longtable}
}

\subsection{Notation}
The notation in the following table is used without defining it by name.
\begin{center}
{\tiny \begin{tabular}
{|l|l|} \hline ${\textbf{C}}$                 & complex numbers \\
\hline ${\textbf{R}}$                 & real numbers \\
\hline ${\textbf{Z}}$                 & integers \\
\hline ${\textbf{F}}$                 & numberfield \\
\hline ${\textbf{P}}^n$               & complex projective space \\
\hline $\textrm{Pic}({\textrm{X}})$    & the Picard group \\
\hline $|D|$                   & linear series \\
\hline $\varphi_D$             & the associated map of a divisor \\
\hline $p_g({\textrm{X}})$             & the geometric genus \\
\hline $p_a({\textrm{X}})$             & the arithmetic genus \\
\hline $h^i(D)$                & the $i$-th Betti number \\
\hline $\delta_p({\textrm{X}})$        & the delta invariant \\
\hline $V(F(x))$               & the zeroset of some polynomial $F(x)$ \\
\hline $[n]$                   & $\{1,2,\ldots,n\}$ with $n\in{\textbf{Z}}_{>0}$ \\
\hline $A_n, D_n, E_n$         & ADE-singularities or Dynkin type\\
\hline
\end{tabular}
 }
\end{center}
For some often used theorems we use the following abbreviations:
\\
{\tiny \begin{tabular}
{|l|l|} \hline (BZ) & (Bezout's theorem)          \\
\hline (HW) & (Hurwitz formula)           \\
\hline (CT) & (Cliffords theorem)         \\
\hline (GF) & (genus formula)             \\
\hline
\end{tabular}
 }
\\
See chapter 2, section 5 in \cite{nls2} for the exact statement of these theorems.

 The assumptions in a claim are always implicitly assumed in the subclaims of this claim. For example if the claim is: ``Claim[1]: If A then B.'' Then in the subclaims of this claim, A is implicitly assumed.

\section{Weak Del Pezzo surfaces of degree one}
In this section we recall some theory concerning degree one Del Pezzo surfaces. There is a bit of overlap with the introduction, but here we will be more detailed. We summarize the necessary concepts needed in this paper. See \cite{dol1} or \cite{man1} for more info about weak Del Pezzo surfaces and further references. In appendix E in \cite{nls2} also some theory concerning Del Pezzo surfaces is recalled in the same notation as here.

 A \textit{weak Del Pezzo surface}\index{weak Del Pezzo surface} is defined as a complex non-singular surface with nef and big anticanonical class $-K$. The \textit{degree}\index{weak Del Pezzo surface!degree}\index{degree} of a weak Del Pezzo surface is defined as $K^2$.

 We can also define a weak Del Pezzo surface as the projective plane blown up in at most $8$ points  such that  no $4$ points lie on a line and no $7$ points lie on a conic. The pull back of the exceptional curves define the \textit{standard Del Pezzo basis}\index{standard Del Pezzo basis} for the Picard group: $\textrm{Pic}{\textrm{X}}\cong{\textbf{Z}}\langle H,Q_1,\ldots, Q_r\rangle $, where $r=9-K^2$. For $r>1$ we have the following intersection product: $H^2=1$, $Q_iQ_j=-\delta_{ij}$, $HQ_i=0$ for $i\in [r]$.

 We define the \textit{minus two classes}\index{weak Del Pezzo surface!minus two classes}\index{minus two classes} to be the set of classes $C$ in $\textrm{Pic}{\textrm{X}}$  such that  $C^2=-2$ and $CK=0$. For example if three points lie on a line then $H-Q_1-Q_2-Q_3$ is an effective minus two class. Here $Q_1,Q_2$ and $Q_3$ are the pull back of the exceptional curves, which blow down to these three points on a line. Note that indeed $(H-Q_1-Q_2-Q_3)^2=-2$ and $(H-Q_1-Q_2-Q_3)K=0$. In particular we see that $-K$ is not ample. The \textit{minus one classes}\index{weak Del Pezzo surface!minus one classes}\index{minus one classes} are defined as the classes  such that  $C^2=CK=-1$. The $Q_i$ are examples of minus one classes. We can consider minus one and effective minus two classes also as divisors. The reason is that these classes have a unique effective representative.

 For degree one weak Del Pezzo surfaces only a multiple of the anticanonical class defines a birational morphism. The anticanonical model ${\textrm{Y}}$ of a weak degree one Del Pezzo surface ${\textrm{X}}$ is defined as $\ensuremath{{\textrm{X}}\stackrel{\varphi_{-3K}}{\rightarrow}{\textrm{Y}}\subset{\textbf{P}}(1:1:2:3)}$ where $K$ is the canonical divisor class of ${\textrm{X}}$. The \textit{canonical projection of a weak degree one Del Pezzo surface} is defined by the 2:1 map associated to twice anticanonical class $\ensuremath{{\textrm{Y}}\stackrel{\varphi_{-2K}}{\rightarrow}{\textbf{P}}(1:1:2)}$. Note that ${\textbf{P}}(1:1:2)$ is isomorphic to the quadric cone ${\textrm{Q}}\subset {\textbf{P}}^3$. The branching curve of $\varphi_{-2K}$ is a cone curve.

 Hyperplane sections of the anticanonical model ${\textrm{Y}}$ of a weak Del Pezzo surface are projected via $\ensuremath{{\textrm{Y}}\stackrel{\varphi_{-2K}}{\rightarrow}Q\subset{\textbf{P}}^3}$ to hyperplanes which intersect the branching curve. The tritangent plane sections which are tangent to the quadric cone pull back to elliptic curves. The remaining hyperplane sections which are bitangent to the cone curve pull back to rational curves.

 Let $m_p$ be the local intersection multiplicity of a hyperplane $L$ with the branching curve at a point $p$. Let $q$ be the preimage of $p$ under the canonical projection (there is a single point since $p$ lies on the branching curve). Let $H$ be the hyperplane section corresponding to the pullback of $L$ along the canonical projection. We have the following formula for the delta-invariant: $\delta_q(H)=\lfloor \frac{m_p}{2}\rfloor $.

 We have that $H$ is contained in the linear series $|-2K|$. The generic element of $|-2K|$ is a curve of genus two. However if $H$ is bitangent to the cone curve, then $m_p\geq 2$ and we have $\delta_q(H)\geq 1$. It follows that tangent planes of the branching curve pull back to curves with a singularity along the ramification curve. Using the same argument we see that bitangent planes pull back to rational curves with two singularities along the ramification curve. These curves form minimal non-fibration families. Details can be found in \cite{nls2}.

 A singularity on the anticanonical model of a degree one weak Del Pezzo surface lies on the ramification curve and is projected to a singularity on the branch curve of the same Dynkin type.

 This discussion is summarized in the following propositions.
\begin{definition}
\label{def:f2_DP1_ring}
\textrm{\textbf{(DP1 ring)}}
  Let ${\textrm{X}}$ be a degree one weak Del Pezzo surface.  Let $D=-K$ be the anticanonical divisor class of ${\textrm{X}}$.  Let ${\textit{H}}^*$ be the sheaf cohomology functor. A \textit{DP1 ring}\index{DP1 ring} for ${\textrm{X}}$ is defined as a graded algebra $A$  such that  \[ A\cong \overset{}{\underset{i>0}{\oplus}}{\textit{H}}^0({\textrm{X}},iD) \] and
\begin{itemize}\addtolength{\itemsep}{1pt}
\item[$\bullet$] $A={\textbf{C}}[y_0,y_1,y_2,y_3]/\langle F\rangle $  such that  $(y_0,y_1,y_2,y_3)$ has weight $(1,1,2,3)$,
\item[$\bullet$] $F = y_3^2 + G(y_0,y_1,y_2)$ and
\item[$\bullet$] $G = y_2^3 + G_4(y_0,y_1)y_2 + G_6(y_0,y_1)$ is a squarefree form ($\deg G_i=i$).
\end{itemize}

\end{definition}

\begin{proposition}\label{prop:f2_DP1_ring}\textbf{\textrm{(properties of DP1 ring)}}
  Let ${\textrm{X}}$ be a degree one weak Del Pezzo surface.  Let $D=-K$ be the anticanonical divisor class of ${\textrm{X}}$.

\textbf{a)} There exists a DP1 ring for ${\textrm{X}}$.

Let ${\textit{H}}^*$ be the sheaf cohomology functor.
Let $A={\textbf{C}}[y_0,y_1,y_2,y_3]/\langle F\rangle $ be a DP1 ring for ${\textrm{X}}$.

\textbf{b)} We have that

${\textit{H}}^0( {\textrm{X}}, D  ) = {\textbf{C}}\langle ~y_0,~y_1~\rangle $,

${\textit{H}}^0( {\textrm{X}}, 2D ) = {\textbf{C}}\langle ~y_0^2,~y_0y_1,~y_1^2,~y_2~\rangle $, and

${\textit{H}}^0( {\textrm{X}}, 3D ) = {\textbf{C}}\langle ~y_0^3,~y_0^2y_1,~y_0y_1^2,~y_1^3,~y_0y_2,~y_1y_2,~y_3~\rangle $.

\begin{proof}
 See \cite{zar1}.
 \end{proof}
\end{proposition}

\begin{proposition}\label{prop:f2_dp1_D}\textbf{\textrm{(canonical projection of degree one weak Del Pezzo surface)}}
  Let ${\textrm{X}}$ be a degree one weak Del Pezzo surface.  Let $D=-K$ be the anticanonical divisor class of ${\textrm{X}}$.
 Let ${\textrm{C}}$ be the branching curve of $\ensuremath{{\textrm{X}}\stackrel{\varphi_{2D}}{\rightarrow}{\textbf{P}}(1:1:2)}$.

\textbf{a)} We have that $\varphi_{2D}$ is a 2:1 morphism and ${\textrm{C}}$ is a cone curve.

 Let $F_{\geq0}({\textrm{X}})$ be the set of effective minus two classes (seen as divisors).

\textbf{b)} We have that $r$ is an isolated double point on ${\textrm{C}}$  if and only if  $\varphi_{2D}(\cup_i F_i)=r$ for some minus two curves $F_i$ in $F_{\geq0}({\textrm{X}})$.

 Let $L$ be a hyperplane section of ${\textbf{P}}(1:1:2)$ (note that ${\textbf{P}}(1:1:2)$ is isomorphic to the quadric cone in ${\textbf{P}}^3$).
 Let $(L\cdot {\textrm{C}})_r$ denote the local intersection multiplicity of ${\textrm{C}}$ and $L$ at a point $r$.
 Let $G+\ensuremath{\overset{}{\underset{j}{\sum}}}F_j=\varphi_{2D}^*L$ in $|2D|$  such that  $G$ has no minus two curves as components.

\textbf{c)} We have that \[ \delta_p(G)=\lfloor \frac{(L\cdot {\textrm{C}})_r}{2}\rfloor  \] where $r=\varphi_{2D}(p)$.

 Let $E({\textrm{X}})$ be the set of minus one classes seen as divisors (thus the exceptional curves).
 Let $m=(m_i)_i$ be a tuple of nonzero intersection multiplicities of $L$ with ${\textrm{C}}$ , or $m=\infty$ if $L$ is a component of ${\textrm{C}}$.

\textbf{d)} We have that $m$ is in $\{$ $(6)$, $(5, 1)$, $(4, 2)$, $(4, 1, 1)$, $(3, 3)$, $(3, 2, 1)$, $(3, 1, 1, 1)$, $(2, 2, 2)$, $(2, 2, 1, 1)$, $(2, 1, 1, 1, 1)$, $(1, 1, 1, 1, 1, 1)$, $\infty$ $\}$ and either:
\begin{itemize}\addtolength{\itemsep}{1pt}
\item[$\bullet$] $G$ is irreducible and $(0:0:1)\notin L$,
\item[$\bullet$] $G = D_1 + D_2$, $L$ consists of two lines through the vertex $(0:0:1)$,
\item[$\bullet$] $G = 2D_1$, $L$ is a double line through the vertex $(0:0:1)$,
\item[$\bullet$] $G = E + E'$, and $m_i$ is even for all $i$, or
\item[$\bullet$] $G = 2E$ and $L \subset {\textrm{C}}$ is a conic component,
\end{itemize} for some $E,E'\in E({\textrm{X}})$ and $D_1,D_2 \in |D|$.

\begin{proof}
 See appendix E, section 6, proposition 241 in \cite{nls2}.
 \end{proof}
\end{proposition}
From the above proposition we find that the effective minus two classes in the Picard group of a weak degree one Del Pezzo ${\textrm{X}}$ contract to singularities on the branching curve, which is a cone curve. In the next section we discuss the classification of such these singularities.

\section{C1 classification}

\label{sec:f2_C1}The effective minus two classes in the Picard group of a weak degree one Del Pezzo surface form a root subsystem of a root system of Dynkin type $E_8$. We consider two root subsystems $S$ and $S'$ of a root system $R$ to be Weyl equivalent, if there exists an element of the Weyl group of $R$ which sends $S$ to $S'$. There are root subsystems which are not Weyl equivalent but have the same Dynkin type. See \cite{dol1} and \cite{man1} for more information.

 The Weyl equivalence classes of root subsystems of $E_8$ will be represented by a set of minus two classes. Each minus two class will be represented by a \textit{C1 label element}\index{C1 label!C1 label element}\index{C1 label element}. We will define these elements by the means of examples. The more formal definition can be found in \cite{nls2}.

 Let $\textrm{Pic}{\textrm{X}}\cong{\textbf{Z}}\langle H,Q_1,\ldots, Q_8\rangle $ be the standard Del Pezzo basis for the Picard group of a weak degree one Del Pezzo surface ${\textrm{X}}$. The C1 label element of
\begin{itemize}\addtolength{\itemsep}{1pt}
\item[$\bullet$] $Q_1-Q_2$ is $12$,
\item[$\bullet$] $H-Q_1-Q_2-Q_3$ is $1123$,
\item[$\bullet$] $2H-Q_1-Q_2-Q_3-Q_4-Q_5-Q_6$ is $278$, where $7$ are $8$ the indices of the omitted $Q_i$, and
\item[$\bullet$] $3H-2Q_1-Q_2-Q_3-Q_4-Q_5-Q_6-Q_7-Q_8$ is $301$ where $1$ is the index of the $Q_i$ which has coefficient two.
\end{itemize} Up to permutation of the $Q_i$ in the standard Del Pezzo basis we represented all possible minus two classes by C1 label elements

 A \textit{C1 label}\index{C1 label} is defined as a pair $(L,r)$ where $r$ is the rank of the root system (in our case $r=8$) and $L$ a set of C1 label elements as described above. The \textit{C1 classification of $E_8$}\index{C1 classification of $E_8$} is defined as the the set of C1 labels defined by the third column in Table~\ref{tab:f2_coc_cls} (which will be stated in a later section). So each C1 label in this column is a unique representive for the Weyl equivalence class of a root subsystem of $E_8$. In \cite{nls-f1} it is explained how this classification is computed.

 We shall see that not all C1 labels can be realized as a set of effective minus two classes of a weak Del Pezzo surface.

\section{Properties of cone curves}
We recall the definition of a cone curve more formally. We show that a component of a cone curve has geometric genus at most four. We state some properties of divisor classes on cone curves. These properties will be used in later sections to determine the irreducible components of bitangent correspondences of cone curves.

 We refer to the representation of a cone curve in a weighted projective plane as a \textit{weighted cone curve}.

 A cone curve can be related to a weak degree one Del Pezzo surface, which has this cone curve as branching curve. We call this relation the \textit{cone curve to DP1 function}. The singularity configuration is preserved by this relation. It follows that singularities of cone curves are classified by the C1 classification of $E_8$.
\begin{definition}
\label{def:f2_cone_curve}
\textrm{\textbf{(cone curve)}}
 We call ${\textrm{C}}\subset{\textbf{P}}^3$ a \textit{cone curve}\index{cone curve}  if and only if
\begin{itemize}\addtolength{\itemsep}{1pt}
\item[$\bullet$] $\deg{\textrm{C}}=6$,
\item[$\bullet$] ${\textrm{C}}$ is the complete intersection of the quadric cone ${\textrm{Q}}$ and a cubic surface ${\textrm{U}}$,
\item[$\bullet$] ${\textrm{C}}$ is reduced, and
\item[$\bullet$] ${\textrm{C}}$ does not go through the vertex of the cone.
\end{itemize}

\end{definition}

\begin{proposition}\label{prop:f2_cone_curve_div_class}\textbf{\textrm{(properties of divisor classes of cone curve)}}

 Let ${\textrm{C}}\subset{\textbf{P}}^3$ be a cone curve on the quadric cone ${\textrm{Q}}$.
 Let $(p_i)_{i\in[r]}$ on ${\textrm{C}}$ be singular points (possibly infinitely near) with multiplicities $(m_i)_{i\in[r]}$.
 Let $\ensuremath{\tilde{{\textrm{Q}}}\stackrel{\mu}{\rightarrow}{\textrm{Q}}}$ be the resolution of the vertex of ${\textrm{Q}}$ and the singularities of ${\textrm{C}}$.
 Let $\tilde{{\textrm{C}}}$ be the strict transform along $\mu$ of ${\textrm{C}}$.
 Let $H$ be the strict transform along $\mu$ of the hyperplane sections of ${\textrm{Q}}$.
 Let $Q_i$ in $\textrm{Pic}\tilde{{\textrm{Q}}}$  such that  $\mu(Q_i)=p_i$ for $i\in[r]$.
 Let $K_{\tilde{{\textrm{C}}}}$ be the canonical divisor class on $\tilde{{\textrm{C}}}$.

\textbf{a)} We have that $K_{\tilde{{\textrm{C}}}}=(H-\ensuremath{\overset{r}{\underset{i=1}{\sum}}}(m_i-1)Q_i)\cdot\tilde{{\textrm{C}}}$ and $r\leq 4$.

\textbf{b)} We have that $p_a({\textrm{C}})=4$.

\textbf{c)} We have the following table
\begin{center}
{\tiny \begin{tabular}
{|c|c|c|c|c|c|} \hline $p_g({\textrm{C}})$                 & $4$ & $3$ & $2$ & $1$ & $0$\\
\hline number of $g_2^1$ on ${\textrm{C}}$ & $0$ & $0$ & $1$ & $\infty$ & $\infty$ \\
\hline
\end{tabular}
 }
\end{center} where a $g^1_2$ is a (not necessarily complete) linear series of projective dimension one and degree two. Moreover if $p_g({\textrm{C}})=2$ then $K_{{\textrm{C}}}=g^1_2$.

\begin{proof}
 See \cite{nls2}, chapter 6, section 2, subsection 1, proposition 67, page 98.
 \end{proof}
\end{proposition}

\begin{definition}
\label{def:f2_weighted_cone_curve}
\textrm{\textbf{(weighted cone curve)}}
  Let ${\textrm{Q}}: u^2-tv=0$ be the quadric cone.  Let ${\textrm{C}} \subset {\textbf{P}}^3(s:t:u:v)$ be a cone curve on ${\textrm{Q}}$. Let \[ \ensuremath{\mu':{\textbf{P}}(2:1:1)\rightarrow {\textrm{Q}},\quad  (x:y:z)\mapsto (x:y^2:yz:z^2)=(s:t:u:v)} \] be a 2-uple embedding. The \textit{weighted cone curve of ${\textrm{C}}$}\index{weighted cone curve of ${\textrm{C}}$} is defined as ${\textrm{W}}=\mu^{-1}({\textrm{C}}) \subset {\textbf{P}}(2:1:1)$. The \textit{cone curve isomorphism}\index{cone curve isomorphism} $\ensuremath{{\textrm{W}}\stackrel{\mu}{\rightarrow}{\textrm{C}}}$ is defined as $\mu=\mu'|{\textrm{W}}$.

\end{definition}

\begin{proposition}\label{prop:f2_weighted_cone_curve}\textbf{\textrm{(properties of weighted cone curve)}}

 Let $\ensuremath{{\textrm{W}}\stackrel{\mu}{\rightarrow}{\textrm{C}}}$ be the cone curve isomorphism.

\textbf{a)} We have that $\mu$ is an isomorphism.

 Let ${\textrm{W}}: W(x,y,z)=0$.

\textbf{b)} We have that $W$ has weighted degree $6$ and the coefficient of $x^3$ in $W(x,y,z)$ is nonzero.

\begin{proof}
 Left to the reader. Note that \textbf{b)} follows from ${\textrm{C}}$ not going through the vertex of the cone.
 \end{proof}
\end{proposition}

\begin{definition}
\label{def:f2_cone_curve_dp_function}
\textrm{\textbf{(cone curve to DP1 function)}}
  Let $A$ be the set of projective isomorphism classes of cone curves.  Let $B$ be the set of projective isomorphism classes of Del Pezzo pairs of degree one. The \textit{cone curve to DP1 function}\index{cone curve to DP1 function} is defined as \[ \ensuremath{\kappa:A\rightarrow B,\quad [{\textrm{C}}]\mapsto [{\textrm{X}}]} \] where ${\textrm{W}}: W(x,y,z)=0 \subset {\textbf{P}}(2:1:1)$ is the weighted cone curve of ${\textrm{C}}$, and ${\textrm{X}}$ is the degree one weak Del Pezzo surface associated to ${\textbf{C}}[y_0,y_1,y_2,y_3]/\langle y_3^2+W(y_2,y_0,y_1)\rangle $ where $(y_0,y_1,y_2,y_3)$ has weight $(1,1,2,3)$.

\end{definition}

\begin{proposition}\label{prop:f2_cone_curve_dp_function}\textbf{\textrm{(properties of cone curve to DP1 function)}}

 Let $\ensuremath{A\stackrel{\kappa}{\rightarrow}B}$ be the cone curve to DP1 function.

\textbf{a)} We have that $\kappa$ is a well defined isomorphism.

 Let ${\textrm{X}}$ be a degree one weak Del Pezzo surface.
 Let ${\textrm{C}}$ be the branching curve of $\ensuremath{{\textrm{X}}\stackrel{\varphi_{2D}}{\rightarrow}{\textrm{Q}}}$.

\textbf{b)} We have that $[{\textrm{C}}]=\kappa^{-1}([{\textrm{X}}])$.

\textbf{c)} We have that Dynkin diagram of the ADE singularities is preserved by $\kappa$.

\begin{proof}

\textit{Claim:} We have that \textbf{a)}.
\\  Left to the reader.

\textit{Claim:} We have that \textbf{b)}.
\\  It follows from Proposition~\ref{prop:f2_dp1_D}.

\textit{Claim:} We have that \textbf{c)}.
\\  From Proposition~\ref{prop:f2_dp1_D}.b) it follows that all minus two curves on ${\textrm{X}}$ are contracted by $\varphi_{2D}$ to singularities on ${\textrm{C}}$. From Definition~\ref{def:f2_DP1_ring} it follows that after a linear automorphism a fixed singularity is at the origin of ${\textrm{C}}: G(x,y,z)=0$. We have that $y_3^2+G(y_0,y_1,y_2)$ by adding the term $y_3^2$. From Proposition~\ref{prop:f2_ADE_sing} (forward reference) it follows that the Dynkin type of the singularity at the origin is invariant under adding powers of two of coordinate functions.
\end{proof}
\end{proposition}

\section{Components and singularities of cone curves}
We show that a reducible cone curve can only have three conic components or a conic and quartic component. We show that the C1 label of a cone curve determines the degree and geometric genus of each of the components of the cone curve.
\begin{proposition}\label{prop:f2_cone_curve_sing}\textbf{\textrm{(properties of singularities of cone curve)}}

\textbf{a)} We have that the singularity configurations of cone curves are contained by the C1 classification for $E_8$.

 See \textsection\ref{sec:f2_C1} for the C1 classification of $E_8$.

\begin{proof}
 Let ${\textrm{X}}$ be a weak degree one Del Pezzo surface.

\textit{Claim:} We have that \textbf{a)}.
\\ From Proposition~\ref{prop:f2_cone_curve_dp_function} it follows that the effective minus-two classes of ${\textrm{X}}$, determine the Dynkin type of singularities on the branching curve. From Proposition~\ref{prop:f2_dp1_D} we know that the branching curve of ${\textrm{X}}$ is a cone curve. The Weyl equivalence class of the set of effective minus-two classes of ${\textrm{X}}$ is represented by a C1 label in the C1 classification of $E_8$. This claim follows from the discussion in \textsection\ref{sec:f2_C1}.

 See chapter 6, section 2 in \cite{nls2} and \cite{nls-f1} for more details. See also \cite{dol1} for the classification of root systems of weak degree one Del Pezzo surfaces.
 \end{proof}
\end{proposition}

\begin{proposition}\label{prop:f2_cone_curve_comp}\textbf{\textrm{(properties of components of cone curve)}}

 Let ${\textrm{C}} \subset {\textbf{P}}^3$ be a cone curve.

\textbf{a)} We have the following table:
\begin{center}
 \begin{tabular}
{|l|l|l|} \hline $(\deg {\textrm{C}}_i)_{i\in I}$  & $(p_a~{\textrm{C}}_i)_{i\in I}$ & $\ensuremath{\overset{}{\underset{p\in {\textrm{C}}}{\sum}}}\delta_p({\textrm{C}})$         \\
\hline \hline $(2,2,2)$                 & $(0,0,0)$               & $6$         \\
\hline $(4,2)$                   & $(1,0)$                 & $4,5$       \\
\hline $(6)$                     & $(4)$                   & $0,1,2,3,4$ \\
\hline
\end{tabular}

\end{center} where $({\textrm{C}}_i)_{i\in I}$ are the irreducible components of ${\textrm{C}}$.

 Let $(L,8)$ be a C1 label in the C1 classification of $E_8$ (see \textsection\ref{sec:f2_C1}).

\textbf{b)} If $\ensuremath{\overset{}{\underset{p\in {\textrm{C}}}{\sum}}}\delta_p({\textrm{C}})=4$ then $(\deg {\textrm{C}}_i)_{i\in I}=(4,2)$  if and only if  $L$ is in Table~\ref{tab:f2_coc_cls} at index 108, 118, 131, 135  or 153.

\begin{proof}
 Let ${\textrm{W}}: W(x,y,z) \subset{\textbf{P}}(2:1:1)$ be the weighted cone curve of ${\textrm{C}}$.

\textit{Claim 1:} We have that $(\deg {\textrm{C}}_i)_{i\in I}$ is $(6)$, $(2,4)$ or $(2,2,2)$.
\\  From ${\textrm{C}}$ not going through the vertex it follows that there are no lines as components.  From Proposition~\ref{prop:f2_weighted_cone_curve} it follows that $W$ factors as either $(y_2^3+\ldots)$ or $(y_2+\ldots)(y_2^2+\ldots)$ or $(y_2+\ldots)(y_2+\ldots)(y_2+\ldots)$.  It follows that $(3,3)$ is not possible.

\textit{Claim 2:} If either $\deg {\textrm{C}}_i=2,4$ or $6$ then $p_a({\textrm{C}}_i)=0,1$  respectively  $4$.
\\  This claim follows from Proposition~\ref{prop:f2_cone_curve_div_class}.

 Let $N$ be the number of components of the branching curve ${\textrm{C}}$.
 Let $S=\ensuremath{\overset{}{\underset{i\in I}{\sum}}}p_g({\textrm{C}}_i)$ be the sum of the geometric genera.
 Let $T=\ensuremath{\overset{}{\underset{p\in {\textrm{C}}}{\sum}}}\delta_p({\textrm{C}})$ be the sum of delta invariants of the singularities.

\textit{Claim 3:} We have that $T=N-S+3$.
\\  From (GF) it follows that $p_a({\textrm{C}}) - T = S - N + 1$.  From claim 2) it follows that $p_a({\textrm{C}})=4$.

\textit{Claim 4:} We have the following table:
\begin{center}
 \begin{tabular}
{|l|l|l|l|} \hline $(\deg {\textrm{C}}_i)_{i\in I}$  & $N$ & $S$          & $T$ \\
\hline\hline $(2,2,2)$              & $3$ & $0$          & $6$           \\
\hline $(4,2)$                & $2$ & $1,0$        & $4,5$         \\
\hline $(6)$                  & $1$ & $4,3,2,1,0$  & $0,1,2,3,4$   \\
\hline
\end{tabular}

\end{center} for the number of T4 families.
\\  The first two columns follows from claim 1).  The $S$ column follows from claim 2).  The $T$ column follows from claim 3).

\textit{Claim 5:} We have that \textbf{a)}.
\\  This claim follows from claim 4).

 Let $\ensuremath{A\stackrel{\kappa}{\rightarrow}B}$ be the cone curve to DP1 function.
 Let ${\textrm{X}}$ be a representative for $\kappa([{\textrm{C}}])$.
 Let $D=-K$ be the anticanonical divisor class of ${\textrm{X}}$.
 Let $E({\textrm{X}})$ be the set of minus one curves on ${\textrm{X}}$.
 Let $F({\textrm{X}})$ be the set of minus two curves on ${\textrm{X}}$.

\textit{Claim 6:} We have that \textbf{b)}.
\\ From Proposition~\ref{prop:f2_cone_curve_dp_function} it follows that ${\textrm{C}}$ is the branching curve of $\ensuremath{{\textrm{X}}\stackrel{\varphi_{2D}}{\rightarrow}{\textrm{Q}}}$. From Proposition~\ref{prop:f2_dp1_D} it follows that $\varphi_{2D}(E)$ is a line  if and only if  $h^0(D-E)>0$, for all exceptional curves $E$. From Proposition~\ref{prop:f2_dp1_D} it follows that ${\textrm{C}}$ has a conic component  if and only if  $h^0(2D-2E)>0$ and $h^0(D-E)=0$ for some exceptional curves $E$. We computed for each C1 label whether there exists an exceptional curve $E$  such that  $h^0(2D-2E)>0$ and $h^0(D-E)=0$.
\end{proof}
\end{proposition}

\section{Germs}
We recall the classification of simple germs up to right equivalence. We state the delta invariant for ADE singularities. This section is needed to analyze ramifications of projections of cone curves.
\begin{definition}
\label{def:f2_ADE_sing}
\textrm{\textbf{(ADE-singularities)}}
 The \textit{ADE-singularities}\index{ADE-singularities} are defined as follows:
\begin{center}
 \begin{tabular}
{|l|lll|} \hline $A_{k\geq 1}$   & $x_1^{k+1} - x_2^2$    & $+$ & $x_3^2+\ldots+ x_n^2$ \\
 $D_{k\geq 4}$   & $x_1^{k-1} - x_1x_2^2$ & $+$ & $x_3^2+\ldots+ x_n^2$ \\
 $E_6$           & $x_1^3     - x_2^4$    & $+$ & $x_3^2+\ldots+ x_n^2$ \\
 $E_7$           & $x_1^3     - x_1x_2^3$ & $+$ & $x_3^2+\ldots+ x_n^2$ \\
 $E_8$           & $x_1^3     - x_2^5$    & $+$ & $x_3^2+\ldots+ x_n^2$ \\
\hline
\end{tabular}

\end{center} for $n\geq 2$ (if $n=2$ then the remaining terms are omitted).

\end{definition}

\begin{proposition}\label{prop:f2_ADE_sing}\textbf{\textrm{(classification of simple space germs)}}

\textbf{a)} A simple complex germ is right equivalent to exactly one of the ADE-singularities.

\begin{proof}
 See \cite{arn1}.
 \end{proof}
\end{proposition}

\begin{proposition}\label{prop:f2_germ_spacecurve_lin}\textbf{\textrm{(properties of space germs of space curves)}}

 Let ${\textrm{C}}\subset{\textbf{P}}^3$ be an analytic space curve.
 Let $\ensuremath{f:({\textbf{C}},0)\rightarrow ({\textrm{C}}',p'),\quad t\mapsto (a'(t),b'(t),c'(t))}$ be an analytic parametrization of a branch of $({\textrm{C}}',p')$.

\textbf{a)} A simple space germ $({\textrm{C}}',p')$ is linear isomorphic to exactly one of the space germs $({\textrm{C}},p)$ as given in the following table:
\begin{center}
{\tiny \begin{tabular}
{|l|| @{$~~(~$} l @{$~~,~~$} l @{$~~,~~$} l @{$~)~~$}|} \hline $({\textrm{C}},p)$           & $a(t)$          & $b(t)$                           & $c(t)$            \\
\hline\hline non-singular         & $t$             & $t^l                + \Phi(l+1)$ & $t^m + \Phi(m+1)$ \\
\hline $A_k$, $k\geq1$ odd  & $t            $ & $+t^{\frac{k+1}{2}} + \Phi(2)$   & $      \Phi(2)$   \\
                      & $t   + \Phi(2)$ & $-t^{\frac{k+1}{2}} + \Phi(2)$   & $      \Phi(2)$   \\
\hline $A_k$, $k\geq2$ even & $t^2 + \Phi(4)$ & $t^{k+1}            + \Phi(4)$   & $      \Phi(4)$   \\
\hline $D_k$, $k\geq4$ even & $t            $ & $+t^{\frac{k-2}{2}} + \Phi(2)$   & $      \Phi(2)$   \\
                      & $t   + \Phi(2)$ & $-t^{\frac{k-2}{2}} + \Phi(2)$   & $      \Phi(2)$   \\
                      & $      \Phi(2)$ & $t                  + \Phi(2)$   & $      \Phi(2)$   \\
\hline $D_k$, $k\geq5$ odd  & $t^2 + \Phi(4)$ & $t^{k-2}            + \Phi(4)$   & $      \Phi(4)$   \\
                      & $      \Phi(2)$ & $t                           $   & $      \Phi(2)$   \\
\hline $E_6$                & $t^4 + \Phi(6)$ & $t^3                + \Phi(6)$   & $      \Phi(6)$   \\
\hline $E_7$                & $t^3 + \Phi(4)$ & $t^2                + \Phi(4)$   & $      \Phi(4)$   \\
                      & $      \Phi(2)$ & $t                           $   & $      \Phi(2)$   \\
\hline $E_8$                & $t^5 + \Phi(6)$ & $t^3                + \Phi(6)$   & $      \Phi(6)$   \\
\hline
\end{tabular}
 }
\end{center} where
\begin{itemize}\addtolength{\itemsep}{1pt}
\item[$\bullet$] $2\leq l<m$, and
\item[$\bullet$] $\Phi(i)$ denotes a powerseries of order at least $i$.
\end{itemize}

\begin{proof}
 See \cite{nls2}, appendix D, proposition 208, page 227.
 \end{proof}
\end{proposition}

\begin{proposition}\label{prop:f2_ade_delta}\textbf{\textrm{(properties of delta invariant of ADE singularities)}}

 Let ${\textrm{C}}$ be an analytic curve.
 Let $({\textrm{C}},p)$ be a space germ.

\textbf{a)} We have the following table:
\begin{center}
 \begin{tabular}
{|l||c|c|c|} \hline $({\textrm{C}},p)$       &  $A_n$  &  $D_n$  & $E_n$ \\
\hline $\delta_p({\textrm{C}})$ & $\lfloor \frac{n+1}{2}\rfloor $ & $\lfloor \frac{n+2}{2}\rfloor $ &  $\lfloor \frac{n+1}{2}\rfloor $ \\
\hline
\end{tabular}

\end{center} where $\delta_p({\textrm{C}})$ is the delta invariant.

\begin{proof}
 See appendix D, section 2, proposition 210, page 230 in \cite{nls2}.
 \end{proof}
\end{proposition}

\section{Ramifications of projections of space curves}
We analyze in detail the ramifications at singularities of linear projections of space curves. See chapter 2, section 5, page 264 in \cite{jha2} or chapter 4, section 2, page 299 in \cite{har1} for the definition of ramification index.

 These ramifications are used later on to find an upper bound on the number of bitangent plane families of cone curves.

 The projection of a space curve in ${\textbf{P}}^3$ along the tangent line $T$ at point $p$ on the curve to the projective line, can be visualized as planes going through the tangent line. The intersection points of the cone curve with a plane (minus $p$) are projected to a point on the image line. The intersection of the space curve with the osculating plane at $p$ projects $p$ to some point.

 If the plane is tangent to another point on the curve (thus a bitangent plane) then this is a ramification point of the projection map. However the projection map can also have a ramification point if the plane goes through a singular point.

 In order to analyze the ramifications we consider a resolution of the singularities of space curves. We can use this analysis in combination with Hurwitz formula. Hurwitz formula states the number of ramification points of projection maps of non-singular curves in terms of the mapping degree and geometric genus.

 For example a plane intersecting a cusp singularity (Dynkin type is $A_2$), intersects after resolution the curve with multiplicity two at this point. In the case of node singularities (Dynkin type is $A_1$) the branches separate after resolution. Thus in the latter case there is no ramification.

 Note that in the examples above we consider the ramification at a singular point  with respect to  a plane through a generic line, although the ramification occurs at the projection from a tangent line.

 In Theorem~\ref{thm:f2_cone_curve_cor} we consider the ramifications of a projection in order to compute the coincidence points of a correspondence. In Proposition~\ref{prop:f2_cone_curve_dual_surface} we consider the projection from a generic line in order to determine the degree of the tangent developable of (the dual of) a cone curve.
\newpage
\begin{proposition}\label{prop:f2_ram_sng}\textbf{\textrm{(ramification of projections of singularities)}}

 Let ${\textrm{C}}\subset{\textbf{P}}^3$ be an analytic curve.
 Let $({\textrm{C}},p)$ be a singular space germ.
 Let $\ensuremath{{\textrm{C}}\stackrel{\varphi}{\rightarrow}{\textbf{P}}^1}$ be the projection from a generic line ${\textrm{L}}\subset{\textbf{P}}^3$ (generic  with respect to  $({\textrm{C}},p)$).
 Let $\ensuremath{\tilde{{\textrm{C}}}\stackrel{\pi}{\rightarrow}{\textrm{C}}}$ be a local parametrization of the branches (thus a resolution of singularities).

\textbf{a)} The local ramification index for $\varphi\circ\pi$ at $({\textrm{C}},p)$ is classified in the following table:
\begin{center}
 \begin{tabular}
{|c||c|c|c|c|c|c|c|c|} \hline type      & $\emptyset$ & $A_{k}$ & $A_{k+1}$ & $D_{k}$ & $D_{k+1}$ & $E_6$ & $E_7$ & $E_8$ \\
\hline $\deg R$  & $0$         & $0$     & $1$       & $1$     & $0$       & $2$   & $1$   & $2$   \\
\hline
\end{tabular}

\end{center} where
\begin{itemize}\addtolength{\itemsep}{1pt}
\item[$\bullet$] the type row denotes the type of the ADE singularity $({\textrm{C}},p)$ ($\emptyset$ means non-singular),
\item[$\bullet$] $k$ is odd, and
\item[$\bullet$] $\deg R = \ensuremath{\overset{}{\underset{\pi(\tilde{p})=p}{\sum}}} \tau(\tilde{p})$ with  $\tau(\tilde{p})$ is the ramification index of $\varphi\circ\pi$ at $\tilde{p}\in\tilde{{\textrm{C}}}$.
\end{itemize}

\begin{proof}
 Let $({\textrm{C}}_j)_j$ be an indexed set of local branches at $p$.
 Let $\ensuremath{\eta_j:({\textbf{C}},0)\rightarrow ({\textrm{C}}_j,p),\quad t\mapsto (~\eta_{j1}(t)),~\eta_{j2}(t),~\eta_{j3}(t)~)}$ be a local parametrization of the branch ${\textrm{C}}_j$ at $p$ (as in Proposition~\ref{prop:f2_germ_spacecurve_lin}).
 Let $\varphi$ locally at $p$ defined by $(x_1,x_2,x_3)\mapsto c_1x_1+c_2x_2+c_3x_3$ .

\textit{Claim:} We have that claim a).
\\  The ramification index of $\varphi\circ\eta_j$ at $0$ is the order of $c_1\eta_{j1}(t)+c_2\eta_{j2}(t)+c_3\eta_{j3}(t)$ minus one.  This claim follows from Proposition~\ref{prop:f2_germ_spacecurve_lin}.
\end{proof}
\end{proposition}

\section{Correspondences of cone curves}
A bitangent correspondence of a curve can be defined by the family of planes which are bitangent to the curve. We generalize techniques from \cite{jha2}, chapter 2, section 5, to state properties about correspondences of space curves.

 The family of tritangent planes, tangent to the singular cone, defines a component of the bitangent correspondence. The bitangent correspondence without this tritangent component is called the \textit{reduced bitangent correspondence}.

 Cone curves are not generic  with respect to  space curves, since there are infinitely many tritangent planes. However we can still extract some invariants of correspondences of cone curves. The arithmetic genus of the reduced bitangent correspondence of a generic cone curve is 181. We also bound the number of special hyperplanes which can occur in a bitangent family.

 The reduced bitangent correspondence of non-generic cone curves will be discussed in later sections.
\begin{definition}
\label{def:}
\textrm{\textbf{(correspondences of curves)}}
  Let ${\textrm{C}}$ be a non-singular algebraic curve. A \textit{correspondence}\index{correspondence} of ${\textrm{C}}$ is defined as an algebraic curve $F\subset {\textrm{C}}\times{\textrm{C}}$. We denote by $F(a)$ the formal sum of divisors $b$  such that  $(a,b)\in F$.

\end{definition}

\begin{definition}
\label{def:}
\textrm{\textbf{(attributes of correspondences of curves)}}
  Let ${\textrm{C}}$ be a non-singular algebraic curve. Let $F\subset {\textrm{C}}\times{\textrm{C}}$ be a correspondence. The \textit{inverse}\index{correspondence!inverse}\index{inverse} of $F$ is defined as $F'=\{~(p,q)~~|~~(q,p)\in F~\}$. Note that $F'(b)$ the formal sum of divisors $a$  such that  $(a,b)\in F$. The \textit{valency}\index{correspondence!valency}\index{valency} of $F$ is defined as $v(F)\in{\textbf{Z}}$  such that  $|F(p) + v(F)p|=|F(q) + v(F)q|$ for all $p,q\in{\textrm{C}}$, or $\infty$ if no such integer exists. The \textit{valency divisor class}\index{correspondence!valency divisor class}\index{valency divisor class} of $F$ is defined as the divisor class of $\hat{F}:=F(p) + v(F)p$. The \textit{degree}\index{correspondence!degree}\index{degree} of $F$ is defined as $d(F):=\deg F(p)$ for generic $p\in{\textrm{C}}$, or $\infty$ if no unique degree exists. The \textit{type}\index{correspondence!type}\index{type} of $F$ is defined as $(~v(F),~d(F),~d(F')~)$. A \textit{united point}\index{correspondence!united point}\index{united point} of $F$ is defined as a point $p\in{\textrm{C}}$  such that  $p\in F(p)$. A \textit{coincident point}\index{correspondence!coincident point}\index{coincident point} of $F$ is defined as a point $(p,q)\in F$  such that  $q$ occurs in $F(p)$ with multiplicity at least 2.

\end{definition}

\begin{proposition}\label{prop:f2_cor}\textbf{\textrm{(properties of correspondences of curves)}}

 Let ${\textrm{C}}$ be a non-singular algebraic curve.
 Let $F\subset {\textrm{C}}\times {\textrm{C}}$ be a correspondence of type $(v,d,d')\in{\textbf{Z}}^3$.
 Let $p_g({\textrm{C}})=g$.
 Let $\ensuremath{{\textrm{C}}\times {\textrm{C}}\stackrel{\pi_1}{\rightarrow}{\textrm{C}}}$ and $\ensuremath{{\textrm{C}}\times {\textrm{C}}\stackrel{\pi_2}{\rightarrow}{\textrm{C}}}$ be the first and second projection map.
 Let ${\textit{H}}_*$ be the singular homology functor.
 Let $A \in{\textit{H}}_1({\textrm{C}}\times {\textrm{C}})$ be the homology class of a fiber of $\pi_1$.
 Let $B \in{\textit{H}}_1({\textrm{C}}\times {\textrm{C}})$ be the homology class of a fiber of $\pi_2$.
 Let $D \in{\textit{H}}_1({\textrm{C}}\times {\textrm{C}})$ be the homology class of the diagonal in ${\textrm{C}}\times {\textrm{C}}$.
 Let $K \in{\textit{H}}_1({\textrm{C}}\times {\textrm{C}})$ be the homology class of the canonical divisor class on ${\textrm{C}}\times {\textrm{C}}$.

\textbf{a)} We have that
\begin{itemize}\addtolength{\itemsep}{1pt}
\item[$\bullet$] $F \sim_h (d'+v)A+(d+v)B-vD$ and $K \sim_h (2g-2)A+(2g-2)B$.
\end{itemize}

\textbf{b)} We have the following intersection products:
\begin{center}
{\tiny \begin{tabular}
{|l|l|l|l|l|l|} \hline $\cdot$ & $A$ & $B$ & $D$    & $K$        & $F$            \\
\hline \hline $A$     & $0$ & $1$ & $1$    & $2-2g$     & $d$            \\
\hline $B$     &     & $0$ & $1$    & $2-2g$     & $d'$           \\
\hline $D$     &     &     & $2-2g$ & $4-4g$     & $d+d'+2gv$     \\
\hline $K$     &     &     &        & $8(g-1)^2$ & $2(g-1)(d+d')$ \\
\hline $F$     &     &     &        &            & $2(dd'-gv^2)$  \\
\hline
\end{tabular}
 }
\end{center} where the empty entries are determined by symmetry of the intersection product.

\textbf{c)} We have that $p_a(F) = (dd'- gv^2) + (g-1)(d+d') + 1$.

 Let $up(F)$ be the number of united points of $F$ (counted with multiplicity).

\textbf{d)} We have that $up(F)=d+d'+2gv$.

 Let $cp(F)$ be the number of coincident points of $F$.
 Let $D(F)$ be the Dynkin diagram of the singularities of $F$.
 Let $\alpha$ be the number of irreducible components in $D(F)$ of type $A_k$ and $D_{k+1}$ for $k$ odd.
 Let $\beta$ be the number of irreducible components in $D(F)$ of type $E_6$ and $E_8$.

\textbf{e)} We have that $cp(F)=(2p_g(F)-2) - d(2g-2) \alpha - \beta$.

\begin{proof}

\textit{Claim 1:} We have that \textbf{a)} and $DD=2-2g$.
\\  See \cite{jha2}, chapter 2, section 5.

\textit{Claim 2:} We have that \textbf{b)}
\\  It follows from the definitions that $AA=BB=0$ and $DA=DB=AB=1$.  This claim follows from claim 1).

\textit{Claim 3:} We have that \textbf{c)} and \textbf{d)}.
\\  From (GF) it follows that $p_a(F)=\frac{1}{2}FF + \frac{1}{2}FK +1$.  From the definitions it follows that $up(F)=DF$.  This claim follows from claim 2).

 Let $\ensuremath{F\stackrel{f}{\rightarrow}{\textrm{C}}}$ where $f=\pi_1|F$ be the first projection restricted to $F$.
 Let $R$ be the ramification divisor of $f$.
 Let $\ensuremath{\tilde{F}\stackrel{\pi}{\rightarrow}F}$ be a resolution of singularities of $F$.
 Let $\ensuremath{\tilde{F}\stackrel{\tilde{f}}{\rightarrow}{\textrm{C}}}$ where $f=(\pi_1|F)\circ \pi$ be the resolution of singularities followed by the first projection.
 Let $\tilde{R}$ be the ramification divisor of $\tilde{f}$.

\textit{Claim 4:} We have that \textbf{e)}.
\\  From the definition of coincident point it follows that $cp(F)=\deg R$.  We have that $\deg f=\deg \tilde{f}=d$.  From (HW) applied to $\tilde{f}$ it follows that $2p_g(F)-2=d(2g-2)+\deg \tilde{R}$.  We may assume that $\pi_1\circ\pi$ restricted to the singular space germ has the same local ramification as $\varphi_1\circ\pi$ as in Proposition~\ref{prop:f2_ram_sng}.  The singularities which do not contribute to $\tilde{R}$ after the resolution by $\pi$ are also coincident points.  From Proposition~\ref{prop:f2_ram_sng} it follows that there are $\alpha$ such singularities.  The singularities which contribute with multiplicity 2 to $\tilde{R}$ after the resolution by $\pi$ are by definition just one coincident point each (since we do not count with multiplicity).  From Proposition~\ref{prop:f2_ram_sng} it follows that there are $\beta$ such singularities.  It follows that $cp(F)=\deg \tilde{R}+\alpha-\beta$.
\end{proof}
\end{proposition}

\begin{definition}
\label{def:f2_cor_tangent}
\textrm{\textbf{(correspondences of cone curves)}}
 Let ${\textrm{C}}\subset{\textbf{P}}^3$ be a cone curve. The \textit{bitangent correspondence}\index{bitangent correspondence} of ${\textrm{C}}$ is defined as: \[ U=\{~(a,b)~~|~~ \textrm{there exists a plane which is bitangent at } a {\textrm{~and~}} b ~\} \subset {\textrm{C}}\times {\textrm{C}}. \] The \textit{tritangent correspondence}\index{tritangent correspondence} of ${\textrm{C}}$ is defined as: \[ V=\{~(a,b)~~|~~a {\textrm{~and~}} b \textrm{ lie on a tritangent plane of } {\textrm{C}} \textrm{ which is tangent to the cone} ~\} \subset {\textrm{C}}\times {\textrm{C}}. \] The \textit{reduced bitangent correspondence}\index{reduced bitangent correspondence} of ${\textrm{C}}$ is defined as: \[ F=\overline{U-V} \subset {\textrm{C}}\times {\textrm{C}}.  \] The \textit{osculating correspondence}\index{osculating correspondence} of ${\textrm{C}}$ is defined as: \[ G=\{~(a,b)~~|~~ b \textrm{ lies on the 3-osculating plane of } a~\}\subset {\textrm{C}}\times {\textrm{C}}. \]

\end{definition}

\begin{theorem}\label{thm:f2_cone_curve_cor}\textbf{\textrm{(properties of correspondences of cone curves)}}

 Let ${\textrm{C}}\subset{\textbf{P}}^3$ be a generic cone curve.
 Let $U,V,F,G$ be as in Definition~\ref{def:f2_cor_tangent}.

\textbf{a)} We have the following table:
\begin{center}
 \begin{tabular}
{|c||c|c|c|c|} \hline      & $U$  & $V$  & $F$  & $G$     \\
\hline\hline $d$  & $14$ & $ 2$ & $12$ & $ 3$    \\
\hline $d'$ & $14$ & $ 2$ & $12$ & $33$    \\
\hline $v$  & $ 4$ & $ 1$ & $ 3$ & $ 3$    \\
\hline
\end{tabular}

\end{center} where $v(U)$ is the valency of $U$, $d(U)$ is the degree of $U$, and $d'(U)$ is the degree of $U'$. Similar for the columns with correspondences $V$, $F$ and $G$.

 Let $\#(m_1,\ldots,m_l)$ denote the number of hyperplane sections $m_1p_1+\ldots+m_lp_l$ of ${\textrm{C}}$ with $m_i\geq m_{i+1}$ the multiplicities.

\textbf{b)} We have that $p_a(F)=181$.

\textbf{c)} We have that $\#(4,2)=12$, $\#(4,1,1)\leq 48$ and $\#(3,2,1)\leq 288$.

\begin{proof}

\textit{Claim 1:} We have that $U(p)=U'(p)$, $V(p)=V'(p)$ and $F(p)=F'(p)$ for generic $p\in{\textrm{C}}$.
\\  From the definitions it follows that the correspondences are symmetric.  For for non-generic $p$, there could be a point in $U(p)$ with multiplicity of two or higher (so we need to assume that $p$ is generic).

 Let $\hat{U}$ be the valency divisor class of $U$. Similarly for $\hat{V}$, $\hat{F}$ and $\hat{G}$.
 Let ${\textrm{T}}=T_p{\textrm{C}}$ be the tangent line for a generic $p\in{\textrm{C}}$.
 Let $\ensuremath{{\textrm{C}}\stackrel{\varphi}{\rightarrow}{\textbf{P}}^1}$ be the projection map from ${\textrm{T}}$.
 Let $R$ be the ramification divisor of $\varphi$.
 Let $K_{\textrm{C}}$ be the canonical divisor class of ${\textrm{C}}$.
 Let $H$ be the class of hyperplane sections of ${\textrm{C}}$.

\textit{Claim 2:} We have that $v(U)=4$, $d(U)=14$ and $\hat{U}=3H$.
\\  Each ramification point $q\in{\textrm{C}}$ of $\varphi$ corresponds to a bitangent plane of ${\textrm{C}}$ through $p$ and $q$.  It follows that $U(p)=R$.  The fibers of $\varphi$ are sections of hyperplanes through ${\textrm{T}}$ minus $p$ with multiplicity two.  From $\deg{\textrm{C}}=6$ it follows that $\deg\varphi=4$.  From (HW) it follows that $\deg R=14$.  We have that $K_{\textrm{C}}=-2\varphi^*H_{{\textbf{P}}^1}+R$ where $H_{{\textbf{P}}^1}$ are the hyperplane sections of ${\textbf{P}}^1$.  From Proposition~\ref{prop:f2_cone_curve_div_class} it follows that $K_{\textrm{C}}=H$.  It follows that $K_{\textrm{C}}=-2(H-2p)+R$ and thus $U(p)=3H-4p$.

 Let ${\textrm{Q}}$ be the quadric cone on which ${\textrm{C}}$ lies.
 Let $\ensuremath{{\textrm{C}}\stackrel{\varphi_M}{\rightarrow}{\textbf{P}}^2}$ be the projection map from the vertex of ${\textrm{Q}}$ (associated to the divisor class $M$).
 Let $M_p$ be the representative in $M$ containing $p\in{\textrm{C}}$

\textit{Claim 3:} We have that $v(V)=1$, $d(V)=2$ and $\hat{V}=M$.
\\  From $V(p)=M_p-p$ and $M_p=\frac{1}{2}H$ as divisor class it follows that $v(V)=1$.  From (BZ) it follows that the intersection of a hyperplane, which is tangent to ${\textrm{Q}}$ and goes through $p$, equals $2p+2q+2r$ for some $q,r\in{\textrm{C}}$.  From $V(p)=q+r$ it follows that $d(V)=2$.

\textit{Claim 4:} We have that $v(F)=3$, $d(F)=12$ and $\hat{F}=3H-M$.
\\  We have that $F(p)=U(p)-V(p)=(3H-4p)-(M-p)=(3H-M)-3p$ as divisor class.  We have that $F(p)=U(p)-V(p)=R-q-r$  such that  $2p+2q+2r$ is a hyperplane section of ${\textrm{C}}$ as divisor.  It follows that we can subtract valences and degrees of $U(p)$ and $V(p)$.

 Let $H_p'$ be a representative of $H$ which is the osculating plane at $p\in{\textrm{C}}$.

\textit{Claim 5:} We have that $v(G)=3$, $d(G)=3$ and $\hat{G}=H$.
\\  From $G(p)=H_p' - 3p$ it follows that $v(G)=3$.  From (BZ) and $\deg{\textrm{C}}=6$ it follows that $d=\deg G(p)=3$.

 Let $\ensuremath{{\textrm{C}}\subset{\textbf{P}}^3\stackrel{\psi_p}{\rightarrow}{\textrm{A}}\subset{\textbf{P}}^2}$ be the projection from $p\in{\textrm{C}}$.
 Let ${\textrm{A}}^*\subset{\textbf{P}}^{2*}$ be the first associated curve of ${\textrm{A}}$ (see appendix A, section 7 in \cite{nls2} or \cite{jha2}).

\textit{Claim 6:} We have that $d'(G)=33$.
\\  We have that $d'(G)=\deg G'(p)$ is the number of planes through $p$, which are 3-osculating planes of ${\textrm{C}}$ for generic $p\in{\textrm{C}}$.  The flexes of ${\textrm{A}}$ are pulled back to osculating planes of ${\textrm{C}}$ through $p$.  From Pl\"ucker formulas for plane curves in \cite{jha2} or appendix B, section 6, proposition 174 in \cite{nls2}, it follows that $d'(G)$ is equal to the number of cusps of ${\textrm{A}}^*$.  From $\psi_p=\varphi_{H_{{\textbf{P}}^3}-E}|{\textrm{C}}$ after blowing up $p$ and $(H_{{\textbf{P}}^3}-E)^2=H_{{\textbf{P}}^3}^2-1$ it follows that $\deg{\textrm{A}}=5$.  After projecting a non-singular curve only $A_1$ singularities are introduced and thus $b=0$.  From $g({\textrm{A}})=4$ and ${\textrm{C}}$ non-singular it follows that ${\textrm{A}}$ has $2$ nodes.  It follows from Pl\"ucker formulas for plane curves that ${\textrm{A}}^*$ has $33$ cusps.

\textit{Claim 7:} We have \textbf{a)}.
\\  From claim 1) it follows that $d=d'$ for $U,V$ and $F$.  From claim 2), claim 3), claim 4), claim 5) and claim 6) it follows that $v$, $d$ and $d'$ hold.

\textit{Claim 8:} Generically we have that $(m_1,\ldots,m_l)$ equals either $(1,1,1,1,1,1)$, $(2,1,1,1,1)$, $(2,2,1,1)$, $(3,1,1,1)$, $(4,1,1)$, $(3,2,1)$, $(2,2,2)$ or $(4,2)$.
\\  We use the assumption that ${\textrm{C}}$ is generic.  If there are infinitely many hyperplanes for some fixed $(m_1,\ldots,m_l)$ then we assume that there possibly exists (up to permutation) hyperplanes with multiplicities $(m_1,\ldots,m_i+m_j,\ldots,m_l)$ where $m_i$ and $m_j$ are omitted.  Geometrically, this means that the points $p_i$ and $p_j$ may coincide.  There are $\infty$ many $(1,1,1,1,1,1)$ and thus $(2,1,1,1,1)$ is possible.  There are $\infty$ many $(2,1,1,1,1)$ and thus $(2,2,1,1)$ and $(3,1,1,1)$ are possible.  There are $\infty$ many $(2,2,1,1)$ and thus $(2,2,2)$, $(4,1,1)$ and $(3,2,1)$ are possible.  There are $\infty$ many $(2,2,2)$ (since ${\textrm{C}}$ lies on a quadric cone) and thus $(4,2)$ is possible.

 Let $D\subset {\textrm{C}}\times{\textrm{C}}$ be the diagonal correspondence.

\textit{Claim 9:} The following table denotes the points in the correspondences for given multiplicities of hyperplane sections which contribute to the correspondences:
\begin{center}
{\tiny \begin{tabular}
{|l||l|l|l|l|l|l|} \hline         & $(3,1,1,1)$ & $(2,2,1,1)$ & $(4,1,1)$     & $(3,2,1)$     & $(2,2,2)$             & $(4,2)$       \\
\hline\hline $U$,$F$ & $\cdot$     & $(p_1,p_2)$ & $(p_1,p_1)$   & $(p_1,p_2)$   & $(p_1,p_2),(p_2,p_1)$ & $(p_1,p_1)$   \\
         &             & $(p_2,p_1)$ &               & $2(p_2,p_1)$  & $(p_1,p_3),(p_3,p_1)$ & $(p_1,p_2)$   \\
         &             &             &               &               & $(p_2,p_3),(p_3,p_2)$ & $3(p_2,p_1)$  \\
\hline $V$     & $\cdot$     & $\cdot$     & $\cdot$       & $\cdot$       & $``$                  & $``$          \\
\hline $G$     & $(p_1,p_2)$ & $\cdot$     & $(p_1,p_1)$   & $2(p_1,p_2)$  & $\cdot$               & $(p_1,p_1)$   \\
         & $(p_1,p_3)$ &             & $(p_1,p_2)$   & $(p_1,p_3)$   &                       & $2(p_1,p_2)$  \\
         & $(p_1,p_4)$ &             & $(p_1,p_3)$   &               &                       &               \\
\hline
\end{tabular}
 }
\end{center} where $\cdot$ denotes that $m_1p_1+\ldots+m_lp_l$ does not contribute to the correspondence, and $``$ denotes a copy of the above column. For $(p,q)$ in each entry we consider the multiplicity degree of $q$ in $U(p)$, $F(p)$, $V(p)$ and $G(p)$.
\\  We recall that $U(p)=R$, $V(p)=M_p-p$, $F(p)=U(p)-V(p)$ and $G(p)=N_p-3p$.  From claim 8) it follows that all contributing hyperplane sections are considered.  This claim follows from the definitions of the correspondences.  For example if $(4,2)$ then $U(p_2)=3p_1+\ldots$ and thus we denote $3(p_2,p_1)$.

 We remark that from Jonquires formula (see \cite{jha3} or appendix B, section 6, page 202 in \cite{nls2}) it follows that for a generic degree six space curve on a non-singular quadric we have that $\#(4,1,1)=60$, $\#(4,2)=0$, $\#(3,2,1)=324$, $\#(2,2,2)=120$.
 Let $\#_F(m_1,\ldots,m_l)$ denote the number of bitangent hyperplane sections $m_1p_1+\ldots+m_lp_l$ defined by a point in $F$.

\textit{Claim 10:} We have that $\#_F(4,2)=0$ and
\begin{itemize}\addtolength{\itemsep}{1pt}
\item[$\bullet$] $DV=\#(4,2)=12$, $DU=c_0\#(4,1,1)+\#(4,2)=60$, $DF=c_0\#(4,1,1)=48$,
\item[$\bullet$] $FG=c_1\#(4,1,1)+c_2\#(3,2,1) =   360$,
\end{itemize} for $c_0,c_1,c_2\in {\textbf{Q}}_{\geq 1}$.
\\  We match common points of the correspondences as in the table of claim 9).  From $\#(4,2)>0$ being a non-generic phenomenon it follows that intersection multiplicities can be higher than one.  We correct the intersection multiplicities with the coefficients $c_0,c_1,c_2\in {\textbf{Q}}_{\geq 1}$.  We have that $\#(4,2)$ is equal to the ramification degree of $\ensuremath{{\textrm{C}}\stackrel{\varphi_M}{\rightarrow}{\textbf{P}}^2}$.  From (HW) it follows that $\#(4,2)=12$.  From $DV=\#(4,2)$ and claim 9) it follows that the points $(p_1,p_1)\in V$ contributing to $\#(4,2)$ are non-singular and intersect $D$ transversely.  From $V\subset U$ it follows that the points $(p_1,p_1)\in U$ contributing to $\#(4,2)$ also are non-singular and intersect $D$ transversely.  From claim 9) it follows that $DU=c_0\#(4,1,1)+\#(4,2)=60$.  From $F=\overline{U-V}$ it follows that $DF=c_0\#(4,1,1)+c\#_F(4,2)=48$ for some $c\in {\textbf{Q}}_{\geq 1}$.  From $\#(4,2)=12$ it follows that $\#_F(4,2)=0$.

 We note that the points in $F$ which contribute to $FV=c\#_{F\cap V}(2,2,2)= 24$ for some $c\in {\textbf{Q}}_{\geq 1}$ are included by the Zarisky closure of $U-V$.

\textit{Claim 11:} Generically $F$ is irreducible.
\\  Left to the reader.  One could consider the singular locus of the tangent developable of ${\textrm{C}}$ as in Proposition~\ref{prop:f2_cone_curve_dual_surface}.

 Let $cp(F)$ be the number of coincident points of $F$.
 Let $\alpha$ be the number of irreducible components in $D(F)$ of type $A_k$ and $D_{k+1}$ for $k$ odd.
 Let $\beta$ be the number of irreducible components in $D(F)$ of type $E_6$ and $E_8$.

\textit{Claim 12:} We have that $p_g(F)=\frac{1}{2}(\#(3,2,1)+74-\alpha+\beta)$ and $\#(3,2,1)\leq 288-\alpha-3\beta$.
\\  From Proposition~\ref{prop:f2_cor} and claim 11) it follows that $2p_g(F)=cp(F) + 74 - \alpha + \beta$.  From Proposition~\ref{prop:f2_cor} it follows that $p_a(F) = 181$.  From $2p_g(F) \leq 2p_a(F)$ it follows that that $cp(F) + 74 - \alpha + \beta \leq 362 - 2\alpha - 2\beta$.  From claim 9) and claim 10) it follows that $cp(F)=\#(3,2,1)$.

\textit{Claim 13:} We have that \textbf{b)} and \textbf{c)}.
\\  From Proposition~\ref{prop:f2_cor} it follows that $p_a(F) = 181$.  From claim 10) it follows that $\#(4,2)=12$ and  $\#(4,1,1)\leq 48$.  From claim 12) it follows that $\#(3,2,1)\leq 288$.

 We note that from $\#(3,2,1)\leq 288$ it follows that $c_1+c_2>2$ in claim 10).
\end{proof}
\end{theorem}

\section{Tangent developables of cone curves and its duals}
We consider the tangent developable surface ${\textrm{T}}$ of the second associated curve of a cone curve. The second associated cone curve (also called dual cone curve) is defined by the three osculating planes of the cone curve. For associated varieties of curves see \cite{jha2} or \cite{nls2}.

 The points on ${\textrm{T}}$ are dual to tangent planes of the cone curve. A component of the singular locus of ${\textrm{T}}$, is either dual to the osculating, tritangent or bitangent correspondence. If the cone curve has only nodes and cusps, we can determine the degree of the union of components (of the singular locus of ${\textrm{T}}$) dual to the bitangent correspondences, in terms of the number of cusps.

 We also show that there is a set theoretical bijection between the components of the singular locus of the tangent developable (of the cone curve) and ${\textrm{T}}$. We use these results for Algorithm~\ref{alg:f2_analyze_developable_surface} to compute examples of bitangent correspondences (see Example~\ref{ex:f2_dual_tangent_developable} and Example~\ref{ex:f2_analyze_developable_surface}).

 Note that a bitangent plane can be defined by two intersecting tangent lines of the cone curve. A point on the tangent developable of the dual cone curve corresponds to a tangent plane of the cone curve.
\begin{proposition}\label{prop:f2_cone_curve_dual_surface}\textbf{\textrm{(properties of tangent developable of cone curve and its dual)}}

 Let ${\textrm{C}}\subset{\textbf{P}}^3$ be a cone curve.
 Let ${\textrm{X}}\subset{\textbf{P}}^{3*}$ be the second associated curve of ${\textrm{C}}$ (points correspond to 3-osculating planes of ${\textrm{C}}$).
 Let ${\textrm{T}}\subset{\textbf{P}}^{3*}$ be the tangent developable of ${\textrm{X}}$.

\textbf{a)} Points on ${\textrm{T}}$ correspond to tangent planes of ${\textrm{C}}$. A point where $n$ lines tangent to ${\textrm{X}}$ meet is dual to an $n$-tangent plane of ${\textrm{C}}$.

 Let ${\textrm{S}}_i\subset{\textbf{P}}^{3*}$ for $i\in [m]$ be the irreducible components in the singular locus of ${\textrm{T}}$.

\newpage

\textbf{b)} The singular locus of the ${\textrm{T}}$ has (after relabeling) the following components:
\begin{itemize}\addtolength{\itemsep}{1pt}
\item[$\bullet$] ${\textrm{S}}_1$ is the cuspidal curve ${\textrm{X}}$,
\item[$\bullet$] ${\textrm{S}}_2$ is a conic dual to the tritangent planes, and
\item[$\bullet$] ${\textrm{S}}_{i}$ for $i>2$ are dual to bitangent planes (thus define the components of the reduced bitangent correspondence).
\end{itemize}

\textbf{c)} If ${\textrm{C}}$ has only $A_1$ and $A_2$ singularities, then
\begin{itemize}\addtolength{\itemsep}{1pt}
\item[$\bullet$] $\deg{\textrm{S}}_1 = 6g-2\beta+12$,
\item[$\bullet$] $\deg{\textrm{S}}_2 = 2$, and
\item[$\bullet$] $\deg \underset{i>2}{\sqcup}{\textrm{S}}_i = \frac{1}{2}(2g - \beta  + 9)(2g - \beta + 8) - (7g - 2\beta + 18)$.
\end{itemize} where $\beta$ is the number of $A_2$ singularities and $g=p_g({\textrm{C}})$.

\textbf{d)} The tangents to ${\textrm{X}}$ are $n$-secants of ${\textrm{S}}_i$ for some $i>2$  if and only if  the associated component of the reduced bitangent correspondence has $n$ bitangents through each point.

 Let $D({\textrm{C}})$ be the Dynkin diagram of the singularities of ${\textrm{C}}$.
 Let $\beta$ be the number of irreducible components in $D({\textrm{C}})$ of either one of the following types: $A_2$, $A_4$, $A_6$, $A_8$, $D_5$, $D_7$ or $E_7$.
 Let $\gamma$ be the number of irreducible components in $D({\textrm{C}})$ of either one of the following types: $E_6$ or $E_8$.
 Let ${\textrm{Z}}$ be the tangent developable of ${\textrm{C}}$.

\textbf{e)} We have that $\deg{\textrm{Z}}=\deg{\textrm{T}}=2g+10-\beta-2\gamma$.

 Let ${\textrm{M}}_i\subset{\textbf{P}}^{3}$ for $i\in [m']$ be the irreducible components in the singular locus of ${\textrm{Z}}$.

\textbf{f)} We have that $({\textrm{M}}_i)_i$ and $({\textrm{S}}_i)_i$ are bijective.

\begin{proof}
 Let $U_p \subset {\textrm{C}}$ be a small neighborhood around generic $p\in{\textrm{C}}$.
 Let $x\in{\textrm{X}}$ correspond to the three osculating plane at $p$.
 Let $U^*_x \subset {\textrm{X}}$ be the dual neighborhood of $U_p$ around $x$.

\textit{Claim 1:} We have that \textbf{a)}.
\\  Tangent planes at $p$ correspond to planes through the tangent line $T_p{\textrm{C}}$.  Planes through a line correspond to a line in ${\textbf{P}}^{3*}$.  The three osculating plane at $p$ is a special plane through $T_p{\textrm{C}}$.  It follows that ${\textrm{T}}$ is formed by lines through each point of ${\textrm{X}}$.  The intersection of $n$ of these lines correspond to a $n$-tangent plane of ${\textrm{C}}$.  We have to show that these lines are tangent lines of ${\textrm{X}}$.  We have that $U_p$ approximates $T_p{\textrm{C}}$.  We move a plane along $U_p$  such that  the plane is a three osculating plane.  As we move the plane it rotates along $U_p$.  If we make $U_p$ smaller then the plane approximately rotates along $T_p{\textrm{C}}$.  This claim follows from $U^*_x$ approximating $T_x{\textrm{X}}$.

 Let $\ensuremath{{\textrm{C}}\subset{\textbf{P}}^3\stackrel{\rho_q}{\rightarrow}{\textrm{A}}\subset{\textbf{P}}^2}$ be the projection from a generic point $q$ in ${\textbf{P}}^3$ outside ${\textrm{C}}$.
 Let $d=\deg{\textrm{A}}$.
 Let $a$, $b$ be the number of $A_1$  respectively  $A_2$ singularities of ${\textrm{A}}$.
 Let ${\textrm{A}}^*\subset{\textbf{P}}^{*2}$ be the first associated curve of ${\textrm{A}}$ (also called dual curve).
 Similarly we define $a^*$, $b^*$ and $d^*$ for ${\textrm{A}}^*$.

\textit{Claim 2:} We have that $g=p_g({\textrm{A}})$, $d=6$, $b=\beta$ and $a=10-g-\beta$.
\\  From $\rho_q$ being birational it follows that $g=p_g({\textrm{A}})$.  We have that $d=\deg{\textrm{C}}=6$ after projection from outside ${\textrm{C}}$.  We have that $A_2$ singularities are generically projected to $A_2$ singularities.  From \cite{har1}, chapter 4, section 3, it follows that projecting space curves to the plane only introduces $A_1$ singularities.  This claim follows from claim 1) and (GF).

 Let $o,t2,t3$ be the number of  respectively  osculating, bitangent (which are not tritangent), and tritangent planes through $q$.

\textit{Claim 3:} We have that  $o=\deg {\textrm{S}}_1$, $t3=\deg {\textrm{S}}_2$ and $t2=\deg \underset{i>2}{\sqcup}{\textrm{S}}_i$.
\\  The planes through a point in ${\textbf{P}}^3$ correspond to a plane in ${\textbf{P}}^{3*}$.  The degree of a curve in ${\textbf{P}}^{3*}$ is defined by the number of intersections with a generic plane.

\textit{Claim 4:} We have that $\deg {\textrm{S}}_2=2$.
\\  The tritangent family of ${\textrm{C}}$ is defined by the ruling of the quadric cone ${\textrm{Q}}$.  Tritangent planes go through the vertex and thus correspond to a planar curve in ${\textrm{P}}^{3*}$.  The tritangent planes are the osculating planes of a conic, and the first associated curve of a plane conic is again a conic.

\textit{Claim 5:} We have that $\deg {\textrm{S}}_1=b^*$ and $\deg \underset{i>2}{\sqcup}{\textrm{S}}_i = a^*-6$.
\\  The images under $\rho_q$ of 3-osculating planes are flexes of ${\textrm{A}}$.  A plane through $q$ bitangent to ${\textrm{C}}$ is projected by $\rho_q$ to a bitangent line of ${\textrm{A}}$.  We have that $a^*$ and $b^*$ are the number of bitangent lines  respectively  flexes of ${\textrm{A}}$.  Tritangents of ${\textrm{A}}$ are dual to three branches meeting in ${\textrm{A}}^*$ , with delta invariant 3.  The triple point deforms to three $A_1$ singularities.  It follows that each tritangent plane through $q$ adds three to $a^*$.  This claim follows from claim 3) and claim 4).

\textit{Claim 6:} We have that $b^*=6g-2b+12$ and $a^*=\frac{1}{2}(2g + 9- b)(2g + 8 - b) - g - b^*$.
\\  From the Pl\"ucker formulas for plane curves in (see \cite{jha2}) it follows that $d^*=d(d-1)-2a - 3b$ and $g=\frac{1}{2}(d-1)(d-2)-a-b$ (and similar for the dual curve).  From claim 2) it follows that $d^*= 2g + 10 - b$, $a^*=\frac{1}{2}(d^*-1)(d^*-2) - g - b^*$ and $3 b^*=d^*(d^*-1)-6-2 a^*$.

 See also appendix B, section 6, proposition 174 in \cite{nls2} for the Pl\"ucker formulas for plane curves.

\textit{Claim 8:} We have that \textbf{b)} and \textbf{c)}.
\\  The cuspidal curve is dual to the osculating planes of ${\textrm{C}}$.  This claim follows from claim 4), claim 5) and claim 6).

\textit{Claim 9:} We have that \textbf{d)}.
\\  From claim 1) it follows that ${\textrm{S}}_i$ defines an irreducible bitangent correspondences for $i\geq 2$.  The $n$-secants are defined by the lines tangent to ${\textrm{X}}$.

 Let ${\textrm{L}} \subset {\textbf{P}}^3$ be a generic line.
 Let $e=\#({\textrm{L}}\cap {\textrm{Z}})$.
 Let ${\textrm{L}}^* \subset {\textbf{P}}^{3*}$ be the dual (or one associated) line of ${\textrm{L}}$.
 Let $e'=\#({\textrm{L}}^*\cap {\textrm{T}})$.
 Let $\ensuremath{{\textrm{C}}\stackrel{\varphi}{\rightarrow}{\textbf{P}}^1}$ be the projection from ${\textrm{L}}$ (we use the same notation as in Proposition~\ref{prop:f2_ram_sng}).
 Let $\ensuremath{\tilde{{\textrm{C}}}\stackrel{\pi}{\rightarrow}{\textrm{C}}}$ be the resolution of singularities of ${\textrm{C}}$.
 Let $R$ be the ramification divisor of $\varphi\circ\pi$.

\textit{Claim 10:} We have $e=\deg{\textrm{Z}}=\deg{\textrm{T}}$ and $e=\deg R$ minus ramifications coming from singularities.
\\  We have that $\deg{\textrm{Z}}=e$.  We have that $e$ equals the number of tangent lines of ${\textrm{C}}$ intersecting ${\textrm{L}}$.  It follows that $e$ equals the number of planes through ${\textrm{L}}$ which are tangent to ${\textrm{C}}$.  It follows that $e=\deg R$ minus ramifications coming from singularities.  The planes through ${\textrm{L}}$ correspond to a line ${\textrm{L}}^*$ in ${\textbf{P}}^{3*}$.  We have that points in ${\textrm{L}}^* \cap {\textrm{T}}$ are dual to planes through ${\textrm{L}}$ which are tangent to ${\textrm{C}}$.  It follows that $e=e'=\deg{\textrm{T}}$.

\textit{Claim 11:} We have that \textbf{e)}.
\\  From $\deg{\textrm{C}}=6$ it follows that $\deg\eta=6$.  From (HW) it follows that $\deg R = 2g+10$.  From Proposition~\ref{prop:f2_ram_sng} it follows how the singularities ramify after resolution, and we adjust accordingly.  From claim 10) it follows that $e=2g+10-\beta-2\gamma$.

 Let ${\textrm{P}}$ be a plane which is bitangent at $q,r\in{\textrm{C}}$.
 Let ${\textrm{M}}_1$ be the cuspidal curve ${\textrm{C}}$ on ${\textrm{Z}}$.
 Let ${\textrm{M}}_2$ be the component corresponding to the intersections of tritangent lines of ${\textrm{X}}$ (which support the tritangent planes tangent to the cone).

\textit{Claim 12:} We have that \textbf{f)}.
\\  We have that ${\textrm{P}}$ contains $T_q{\textrm{C}}$ and $T_r{\textrm{C}}$.  It follows that $T_q{\textrm{C}}$ and $T_r{\textrm{C}}$ intersect.  It follows that $T_q{\textrm{C}}\cap T_r{\textrm{C}}$ is a singular point of ${\textrm{Z}}$.  We have that the duals (or one associated curves) of $T_q{\textrm{C}}$ and $T_r{\textrm{C}}$ are lines in ${\textbf{P}}^{3*}$, which intersect in a point $s\in{\textbf{P}}^{3*}$.  We have that ${\textrm{P}}$ is dual to $s$ and thus $s$ is the intersection of two tangent lines of ${\textrm{X}}$.  Conversely the plane spanned by these two tangent lines of ${\textrm{X}}$ is dual to $T_q{\textrm{C}}\cap T_r{\textrm{C}}$.  It follows that there is a birational relation between the singular loci of ${\textrm{Z}}$ and ${\textrm{T}}$.  From this relation being algebraic it follows that irreducible components correspond to irreducible components.  In particular ${\textrm{M}}_1$ is send to ${\textrm{S}}_1$ and ${\textrm{M}}_2$ is send to ${\textrm{S}}_2$.
\end{proof}
\end{proposition}

\section{Bitangent correspondences of irreducible cone curves}
We give an upper bound for the number of irreducible components of the reduced bitangent correspondence of an irreducible cone curve in terms of the singularities.

 The idea of the proof is to consider the number of bitangent planes through a generic point on the cone curve. At a generic point, each bitangent plane belongs to at most one family. It could be that one family contains more than one bitangent plane through some fixed point; I was not able to exclude this.

 The planes through the tangent line at a generic point on the cone curve define a linear projection map. The ramification points of this map correspond to a plane which contains the tangent line of another point on the cone curve. However if the other point is singular then this plane is not a valid bitangent plane.

 We use the analysis of ramifications at singularities of the cone curve and Hurwitz formula to determine the number of (valid) bitangent planes through a generic point on a cone curve with given singularities.
\begin{proposition}\label{prop:f2_cone_curve_bound}\textbf{\textrm{(properties of number of components of reduced bitangent correspondence: (6))}}

 Let $S$ be a set of cone curves with a fixed set of ADE singularities.
 Let ${\textrm{C}} \subset {\textbf{P}}^3$ be a generic cone curve in $S$ on quadric cone ${\textrm{Q}}$.
 Let $D({\textrm{C}})$ be the Dynkin diagram of the singularities of ${\textrm{C}}$.
 Let $k$ be the number of irreducible components in $D({\textrm{C}})$ of either one of the following types: $A_2$, $A_4$, $A_6$, $A_8$, $D_5$, $D_7$ or $E_7$.
 Let $l$ be the number of irreducible components in $D({\textrm{C}})$ of either one of the following types: $E_6$ or $E_8$.
 Let $F\subset {\textrm{C}}\times{\textrm{C}}$ be the reduced bitangent correspondence.

\textbf{a)} We have that \[ 12-2\ensuremath{\overset{}{\underset{p\in{\textrm{C}}}{\sum}}}\delta_p({\textrm{C}})-k-2l \]  is an upper bound for the number of one dimensional components of $F$.

\textbf{b)} If the upper bound is not zero then $F$ has at least one component of dimension one.

\begin{proof}
 Let $({\textrm{C}},p)$ be a space germ where $p\in {\textrm{C}}$ is generic.
 Let ${\textrm{T}}=T_p{\textrm{C}}$ be the tangent line.
 Let $m$ be the number of planes through ${\textrm{T}}$ which are bitangents of ${\textrm{C}}$ but not tangent to ${\textrm{Q}}$.
 Let $n$ be the number of components of the reduced bitangent correspondence $F$.
 Let $G=\{~ (q,c) ~~|~~ (q,r)\in F_i\cap F_j \textrm{ where } 0\leq i\neq j \leq n ~\} \subset {\textrm{C}}\times {\textrm{C}}$

\textit{Claim 1:} We have that $n\leq m$.
\\  We have that $G$ is a finite set and $(p,c) \notin G$ for all $c\in{\textrm{C}}$.  It follows that each bitangent plane through ${\textrm{T}}$ belongs to at most one $F_i$.  Each $F_i$ contains an element corresponding to a bitangent plane through $p$.

 Let $\ensuremath{{\textrm{C}}\stackrel{\varphi}{\rightarrow}{\textbf{P}}^1}$ be the projection map from ${\textrm{T}}$.
 Let $R$ be the ramification divisor of $\varphi\circ\pi$.

\textit{Claim 2:} We have that $m\leq \deg R$.
\\  Each bitangent plane through ${\textrm{T}}$ gives rise to ramification of at least one.

\textit{Claim 3:} We have that $m=12-2\ensuremath{\overset{}{\underset{p\in{\textrm{C}}}{\sum}}}\delta_p({\textrm{C}})-k-2l$
\\  From (GF) it follows that $p_g({\textrm{C}})=4-\ensuremath{\overset{}{\underset{p\in{\textrm{C}}}{\sum}}}\delta_p({\textrm{C}})$.  We have that $\deg\varphi=4$.  From (HW) it follows that $\deg R=14-2\ensuremath{\overset{}{\underset{p\in{\textrm{C}}}{\sum}}}\delta_p({\textrm{C}})$.  From claim 2) it follows that $m\leq \deg R$.  The tritangent plane though ${\textrm{T}}$ tangent ${\textrm{Q}}$ accounts for ramification of degree $2$.  From Proposition~\ref{prop:f2_ram_sng} it follows that the ramification of a generic plane through a singularity is $2$ for $E_6$ or $E_8$.  Similarly we have ramification degree $1$ for $A_2$, $A_4$, $A_6$, $A_8$, $D_5$, $D_7$ and $E_7$.  It follows that we have ramification of degree $k+2l$ coming from planes through ${\textrm{T}}$ and a singularity.  These planes are not bitangent planes.

\textit{Claim 5:} We have that \textbf{a)}.
\\  This claim follows from claim 1) and claim 3).

\textit{Claim 6:} We have that \textbf{b)}.
\\  If $m\neq 0$ then there is a bitangent plane through each point of ${\textrm{C}}$.  It follows that these planes must be defined by an algebraic equation which determines a component of $F$.
\end{proof}
\end{proposition}
We give an upper bound for the number of rational components of the reduced bitangent correspondence. We are still considering the case that the cone curve is irreducible. A rational component of the bitangent correspondence defines a $g^1_2$ on the cone curve. It follows that if the cone curve has genus three or four, there are no rational components and if genus two, then there is at most one component.

 For the genus one and zero case we use the tangent developable of the dual curve. We project the cuspidal curve on the developable to the remaining components of the singular locus, along the ruling. Then we use Hurwitz theorem to say something about the genus of the components if we know the number of components. Recall that the rational components of the reduced bitangent correspondence define unirational parametrizations of minimal lexicographic degree.

 In particular we consider cone curves of genus two. We show that a C1 label (denoting the singularities) is for these curves not sufficient to determine whether it has a rational bitangent correspondence, and give some conditions on the coordinates.
\begin{proposition}\label{prop:f2_cone_curve_correspondence_genus}\textbf{\textrm{(properties of genus of components of reduced bitangent correspondence)}}

 Let ${\textrm{C}} \subset {\textbf{P}}^3$ be a cone curve.
 Let $F \subset {\textrm{C}}\times {\textrm{C}}$ be the reduced bitangent correspondence.
 Let $u$ be the upper bound for the number of components of $F$ as given in Proposition~\ref{prop:f2_cone_curve_bound}.

\textbf{a)} We have the following table:
\begin{center}
{\tiny \begin{tabular}
{|c|c|c|c|c|c|} \hline $p_g({\textrm{C}})$ & $4$ & $3$ & $2$ & $1$ & $0$\\
\hline $r$         & $0$ & $0$ & $1$ & $u$ & $u$ \\
\hline
\end{tabular}
 }
\end{center} where $r$ is an upper bound for the number of rational components.

\textbf{b)} If $n=u$ then the genus of each of the components is bounded by the geometric genus $p_g({\textrm{C}})$.

\begin{proof}
 Let $A={\textrm{C}}\otimes{\textrm{C}}$ be the symmetric tensor of ${\textrm{C}}$ with itself (or in other words the divisors in $\textrm{Pic}({\textrm{C}})$ of degree two).
 Let $F=F_1+\ldots+F_n$ be the decomposition of $F$ into $n$ irreducible components.

\textit{Claim 1:} We have that \textbf{a)}.
\\  Suppose that $F_i$ is rational for some $i\in[n]$.  We have that $F_i\subset A$ is a rational curve on the surface $A$.  It follows that $F_i$ corresponds to a family of degree two divisors over the projective line.  Effective divisors are linear equivalent if they are contained in a family over the projective line.  It follows that $F_i$ defines a $g_2^1$.  From Proposition~\ref{prop:f2_cone_curve_div_class} it follows that this claim holds.

 Let ${\textrm{X}}\subset{\textbf{P}}^{3*}$ (second associated curve of ${\textrm{C}}$), ${\textrm{T}}\subset{\textbf{P}}^{3*}$ (tangent developable of ${\textrm{X}}$) and $S_i \subset{\textbf{P}}^{3*}$ (singular locus of ${\textrm{T}}$) be defined as in Proposition~\ref{prop:f2_cone_curve_dual_surface}.

\textit{Claim 2:} We have that \textbf{b)}.
\\  If $u=n$ then every point on $S_i$ for some $i>2$ is reached by exactly two tangent lines of ${\textrm{X}}$.  It follows that ${\textrm{T}}$ defines a $2:1$ morphism on ${\textrm{S}}_i$ for all $i$.  From (HW) it follows that $p_g(S_i)\leq p_g({\textrm{X}})$.  From ${\textrm{X}}$ and ${\textrm{C}}$ being birational it follows that $p_g({\textrm{X}})=p_g({\textrm{C}})$.
\end{proof}
\end{proposition}

\begin{proposition}\label{prop:f2_cone_curve_genus_two}\textbf{\textrm{(properties of bitangent correspondence of genus two cone curve)}}

 Let ${\textrm{C}}\subset{\textbf{P}}^3$ be a cone curve of geometric genus $2$.
 Let ${\textrm{W}}: W(x,y,z)$ be the weighted cone curve of ${\textrm{C}}$ (see Definition~\ref{def:f2_weighted_cone_curve}).
 Let $F\subset{\textrm{C}}\times{\textrm{C}}$ be the bitangent correspondence.

\textbf{a)} We have that \[ {\textrm{W}}: W(x,y,z)=x^3+P(y,z)x^2+Q(y,z)x+R(y,z) \subset {\textbf{P}}(2:1:1) \] where
\begin{itemize}\addtolength{\itemsep}{1pt}
\item[$\bullet$] $P(y,z)=\ensuremath{\overset{2}{\underset{i=0}{\sum}}}p_i y^{2-i}z^i$,
\item[$\bullet$] $Q(y,z)=\ensuremath{\overset{4}{\underset{i=0}{\sum}}}q_i y^{4-i}z^i$, and
\item[$\bullet$] $R(y,z)=\ensuremath{\overset{6}{\underset{i=0}{\sum}}}r_i y^{6-i}z^i$.
\end{itemize} for $p_i,q_i,r_i$ in ${\textbf{C}}$.

\textbf{b)} If $C$ has singularities $2A_2, 2A_1$ or $A_1+A_2$ then up to projective isomorphism we have that $\{~q_0,~q_4,~r_0,~r_1,~r_5,~r_6~\}$ vanish.

\textbf{c)} If $C$ has singularities $A_3$ or $A_4$ then up to projective isomorphism we have that $\{~q_3,~q_4,~r_3,~r_4,~r_5\}$ vanish.

\textbf{d)} If $C$ has singularities $2A_2, 2A_1$ or $A_1+A_2$ as above then $F$ has one rational component  if and only if  $\{~p_2 q_1 - p_0 q_3,~p_2 r_2 - p_0 r_4,~q_3 r_2 - q_1 r_4~\}$ vanish.

\textbf{e)} If $C$ has singularities $A_3$ or $A_4$ as above then $F$ has one rational component  if and only if
\begin{itemize}\addtolength{\itemsep}{1pt}
\item[$\bullet$] $\{~p_2 q_1 - p_1 q_2,~p_2 r_1 - p_1 r_2,~q_2 r_1 - q_1 r_2~\}$ vanish, or
\item[$\bullet$] $\{~p_2,~ p_1^2 - 4 p_0 p_2 - 4 q_2,~ -4 p_2 q_0 + 2 p_1 q_1 - 4 p_0 q_2 - 4 r_2,~ q_1^2 - 4 q_0 q_2 - 4 p_2 r_0 + 2 p_1 r_1 - 4 p_0 r_2,~ -4 q_2 r_0 + 2 q_1 r_1 - 4 q_0 r_2,~r_1^2 - 4 r_0 r_2~\}$ vanish.
\end{itemize}

\begin{proof}

\textit{Claim 1:} We have that \textbf{a)}.
\\  This claim follows from Proposition~\ref{prop:f2_weighted_cone_curve}.

\textit{Claim 2:} If $p_g({\textrm{C}})=2$ then ${\textrm{C}}$ has $2A_1$, $2A_2$, $A_1+A_2$, $A_3$ or $A_4$ singularities.
\\  This claim follows from (GF) and Proposition~\ref{prop:f2_ade_delta}.

\textit{Claim 3:} We have that \textbf{b)}.
\\  After coordinate change we may assume that $a=(0:0:1)$ and $b=(0:1:0)$ are double points.  This claim follows from computing the Groebner basis of the partial derivatives of $W$ evaluated at $a$ and $b$: ${\textbf{C}}\langle ~ W_x(a) ,~ W_y(a) ,~ W_z(a) ,~ W_x(b) ,~ W_y(b) ,~ W_z(b) ~\rangle .$

 Let $W_0(x,y)=W(x,y,1)$, $W_1(x,y)=\frac{1}{y^2}W_0(xy,y)$, and $W_2(x,y)=\frac{1}{y^2}W_1(xy,y)$.

\textit{Claim 4:} We have that \textbf{c)}.
\\  After coordinate change we may assume that there is a double point at $(0:0:1)$ and an infinitely near point, such that $W_0$, $W_1$, $W_2$ are equations of blowup charts.  See chapter 8, section 1-2 in \cite{nls2} for the description of blowups in terms of charts.  We may assume that $W_1(0,0)=0$ and $W_1,W_2$ are irreducible polynomials.  It follows that this claim holds.

 In the remainder of this proof we assume that ${\textrm{C}}$ has singularities $2A_2, 2A_1$ or $A_1+A_2$ at $c'=(0:0:0:1)$ and $d'=(0:1:0:0)$.
 Let $S=\{~(a,b)\in {\textrm{C}}\times {\textrm{C}}~~|~~ {\textrm{C}}\cap {\textrm{H}}=\{a,b,c,d\} \textrm{ for a hyperplane } {\textrm{H}} \textrm{ through } c' {\textrm{~and~}} d' ~\}^-$.

\textit{Claim 5:} We have that $F$ has a rational component $F_1$  if and only if  $F_1=S$.
\\  We have that $F\subset{\textrm{C}}\times{\textrm{C}}$ is a one dimensional symmetric correspondence of degree 2.  It follows that if $F_1$ is rational then $F_1$ is isomorphic to a $g^2_1$.  The unique $g^2_1$ on ${\textrm{C}}$ is given by the canonical divisor.  From Proposition~\ref{prop:f2_cone_curve_div_class} it follows that $K_{\tilde{C}}=(H-Q_1-Q_2)\cdot\tilde{C}$ where $Q_1$ and $Q_2$ blow down to  respectively  $c'$ and $d'$.

 Let ${\textrm{H}}_l: H_l(x,y,z) = x - l y z=0 \subset{\textbf{P}}(2:1:1)$ for $l\in{\textbf{C}}$.
 Let ${\textrm{T}}_t: T_t(x,y,z) = t_0  x + t_1  y^2 + t_2  y  z + t_3  z^2=0 \subset{\textbf{P}}(2:1:1)$ with coefficients $t=(t_0,\ldots,t_3)\in{\textbf{C}}^4$ .
 Let $c=(0:0:1)$ and $d=(0:1:0)$ (note that $c'=\mu(c)$ and $d'=\mu(d)$).
 Let ${\textrm{W}}\cap {\textrm{H}}_l=\{ a_l,b_l, c,d \}$ for generic $l\in{\textbf{C}}$ (we may assume that $t_0=1$).
 Note that we use the isomorphism $\mu$ to compute in ${\textbf{P}}(2:1:1)$. We see that ${\textrm{H}}_l$ is the pullback along $\mu$ of hyperplanes through $c'$ and $d'$ parametrized by $l\in{\textbf{C}}$ We have that ${\textrm{T}}_t$ is the pullback along $\mu$ of hyperplanes in ${\textbf{P}}^3$, for all $t\in{\textbf{C}}^4$.
 We consider the sections with ${\textrm{W}}$ of a one dimensional linear series of hyperplanes ${\textrm{H}}_l$ through double points $c$ and $d$, which are parametrized by $l$. The hyperplanes intersect ${\textrm{W}}$ in two other points $a_l$ and $b_l$. We want to find conditions on $p_i,q_i$ and $r_i$  such that  the planes ${\textrm{T}}_t$ through $a_l$ and $b_l$ are bitangent planes for all $l$. From claim 5) it follows that if there is a rational component of the bitangent correspondence on ${\textrm{W}}$  if and only if  these conditions on $p_i,q_i$ and $r_i$ are satisfied.
 Let $\alpha_t(p)=\det{\tiny \begin{bmatrix}
 T_x(p) & T_y(p) \\
 W_x(p) & W_y(p) \\

\end{bmatrix}
 } $ for a point $p\in {\textrm{W}}$ and $t\in{\textbf{C}}^4$.
 Let $\ensuremath{f_l:{\textbf{C}}^4\rightarrow {\textbf{C}}^4,\quad t\mapsto  (~ T_t(a_l),~ T_t(b_l),~ \alpha_t(a_l),~ \alpha_t(b_l) ~) }$ be a function depending on $l\in {\textbf{C}}$.
 Let $J(f_l)$ be the Jacobian matrix of $f_l$.
 Let $G_l(x,y,z)=W(lyz,y,z)=G_0(l)y^2 + G_1(l)yz + G_2(l)z^2$ for $l\in{\textbf{C}}$.

\textit{Claim 6:} We have that ${\textrm{T}}_t$ is tangent to ${\textrm{W}}$ at $a_l$ and $b_l$  if and only if  $\det J(f_l)=0$ and $G_l(a_l)=G_l(b_l)=0$.
\\  We have that ${\textrm{T}}_t$ is tangent to ${\textrm{W}}$ at $p$  if and only if  $\alpha_t(p)=0$, $T_t(p)=0$ and $W(p)=0$ for all $t\in{\textbf{C}}^4$.  From $H_l(a_l)=H_l(b_l)=0$ (thus $x=lyz$) and $G_l(a_l)=G_l(b_l)=0$ it follows $W(a_l)=W(b_l)=0$.  We have that $f_l(t)=0$  if and only if  $t \in \ker J(f_l)$.  We have $f_l(t)=0$ for non-trivial $t$  if and only if  $\det J(f_l)=0$.

 Let $Z=G_0(l)( a_{ly} + b_{ly} ) + G_1(l)$ where $a_l=(a_{lx}:a_{ly}:a_{lz})$ and $b_l=(b_{lx}:b_{ly}:b_{lz})$.

\textit{Claim 7:} We have that \textbf{d)}.
\\  We consider the chart $z\neq 0$ and assume $a_{lz}=b_{lz}=1$. We require that $H_l(a_l)=H_l(b_l)=0$. From claim 6) it follows that we require $\det J(f_l)=0$ and $G_l(a_l)=G_l(b_l)=0$. From $G_l(a_l)=G_l(b_l)=0$ it follows that $G_l(x,y,1)=(y-a_{ly})(y-b_{ly})=y^2-(a_{ly} + b_{ly})y+a_{ly}b_{ly}$, and thus $Z=0$. We elimate $a_{lx},a_{ly}$ and $b_{ly}$ using resultants: $R_0=\textrm{res}_{a_{lx}}(\det J(f_l), H_l(a_l))$, $R_1=\textrm{res}_{a_{ly}}(R_0,Z)$ and $R_2=\textrm{res}_{b_{ly}}(R_1,G_l(b))$. We compute the coefficients of $R_2$  with respect to  $l$ with a computer algebra system. This claim follows from claim 5).

\textit{Claim 8:} We have that \textbf{e)}.
\\  The proof is similar as in claim 7).
\end{proof}
\end{proposition}

\section{Bitangent correspondences of reducible cone curves}
We recall that a reducible cone curve consists of either three conic components or a conic and quartic component.

 If there are three conic components then the reduced bitangent correspondence has three rational components. The idea of the proof is to consider two of the conic components of the cone curve. We find that through each tangent line of one of the conics there are two bitangent planes which are tangent to the other conic. One of these planes must be a tritangent plane, tangent to the singular cone. It follows that this gives rise to one of the rational components of the reduced bitangent correspondence. The other two components are obtained in the same way with the two other pairs of conics.
\begin{proposition}\label{prop:f2_cone_curve_222}\textbf{\textrm{(properties of bitangent correspondence: (2,2,2))}}

 Let ${\textrm{C}}\subset{\textbf{P}}^3$ be a cone curve with the three conics ${\textrm{C}}_1,{\textrm{C}}_2$ and ${\textrm{C}}_3$ as components.
 Let $F \subset {\textrm{C}}\times{\textrm{C}}$ be the reduced bitangent correspondence.

\textbf{a)} We have that $F=F_1+F_2+F_3$  such that  $p_g(F_1)=p_g(F_2)=p_g(F_3)=0$.

\begin{proof}
 See \cite{nls2}, chapter 6, section 2, subsection 8, proposition 81, page 130.
 \end{proof}
\end{proposition}
If the cone curve contains a quartic and a conic component, the reduced bitangent component has four components. The genus of one component is equal tho the genus of quartic component of the cone curve, and the other three components are rational. From the singularity configuration of the cone curve it is possible to determine the genus of the quartic component of the cone curve. For the proof we consider the linear series of quadric surfaces, with the quartic component as base locus. We consider the unique quadric in the linear series which contains a stationary bisecant of the quartic component. It follows that this is a quadric cone and it can be shown that there are exactly four cones in this linear series of quadric surfaces. Each cone defines a bitangent family of the quartic component of the cone curve. One of these cones defines a the tritangent family, and thus the quartic has three reduced bitangent families. Moreover these components of the reduced bitangent correspondence are rational. The fourth component of the reduced bitangent correspondence comes from the family of bitangent planes which are bitangent to a point on the conic and a point on the quartic component of the cone curve. The genus of this component is determined by the genus of the quartic component of the cone curve.
\begin{proposition}\label{prop:f2_cone_curve_42}\textbf{\textrm{(properties of bitangent correspondence: (4,2))}}

 Let ${\textrm{C}}\subset{\textbf{P}}^3$ be a cone curve with irreducible components ${\textrm{C}}_1$ and ${\textrm{C}}_2$  such that  $\deg {\textrm{C}}_1=4$ and $\deg {\textrm{C}}_2=2$.
 Let $F \subset {\textrm{C}}\times {\textrm{C}}$ be the reduced bitangent correspondence.

\textbf{a)} We have that $F=F_1+F_2+F_3+F_4$  such that  $p_g(F_i)=0$ for $i\in [4]$ and $p_g(F_4)=p_g({\textrm{C}}_1)$.

\begin{proof}
 See \cite{nls2}, chapter 6, section 2, subsection 8, proposition 82, page 131.
 \end{proof}
\end{proposition}

\section{Classification of bitangent correspondence of cone curves}
We present a table where for each C1 label of a cone curve we give:
\begin{itemize}\addtolength{\itemsep}{1pt}
\item[$\bullet$] the Dynkin type and the sum of delta invariants of the singularities,
\item[$\bullet$] the degree of the components of the cone curve,
\item[$\bullet$] the number of components of the reduced bitangent correspondence and
\item[$\bullet$] the number of rational components of the reduced bitangent correspondence.
\end{itemize} For some of these numbers we can only provide an interval. Two different C1 labels might have the same Dynkin type. Also we note that not all C1 labels define singularities of cone curves. These rows are included with a dash in the columns.

 It is remarkable that the reduced bitangent correspondence of a cone curve with four cusps is the empty set!
\begin{theorem}\label{thm:f2_coc_cls}\textbf{\textrm{(properties of table of components of bitangent correspondence of cone curves)}}

\textbf{a)} We have that Table~\ref{tab:f2_coc_cls} is correct.

\begin{proof}

\textit{Claim 1:} We have that \textbf{a)}.
\\  The `type' and `L' columns follow from Proposition~\ref{prop:f2_cone_curve_sing}. The `delta' column follows from Proposition~\ref{prop:f2_ade_delta}. From Proposition~\ref{prop:f2_cone_curve_comp} it follows that $\ensuremath{\overset{}{\underset{p\in {\textrm{C}}}{\sum}}}\delta_p({\textrm{C}})\leq 6$ and thus no cone curve exists with given singularity configuration, if the `delta' column contains an entry higher than 6. The `deg' column with $(\deg{\textrm{C}}_i)_i$ follows from Proposition~\ref{prop:f2_cone_curve_comp}. The `T5' column follows from Proposition~\ref{prop:f2_cone_curve_bound} if $(\deg{\textrm{C}}_i)_i=(6)$, Proposition~\ref{prop:f2_cone_curve_222} if $(\deg{\textrm{C}}_i)_i=(2,2,2)$ and Proposition~\ref{prop:f2_cone_curve_42} if $(\deg{\textrm{C}}_i)_i=(4,2)$. The `T5R' column follows from Proposition~\ref{prop:f2_cone_curve_correspondence_genus} if $(\deg{\textrm{C}}_i)_i=(6)$, Proposition~\ref{prop:f2_cone_curve_222} if $(\deg{\textrm{C}}_i)_i=(2,2,2)$ and Proposition~\ref{prop:f2_cone_curve_42} if $(\deg{\textrm{C}}_i)_i=(4,2)$.
\end{proof}
\end{theorem}

\begin{tab}\label{tab:f2_coc_cls}\textbf{\textrm{(table of components of bitangent correspondence of cone curves)}}
\begin{itemize}\addtolength{\itemsep}{1pt}
\item[$\bullet$] Let $F\subset{\textrm{C}}\times{\textrm{C}}$ be the reduced bitangent correspondence of a cone curve ${\textrm{C}}\subset{\textbf{P}}^3$.
\item[$\bullet$] The column `type' denotes the Dynkin type of the Weyl equivalence class of the root subsystem in $E_8$.
\item[$\bullet$] The column `L' denotes the C1 label $(L,8)$ for the corresponding root subsystem.
\item[$\bullet$] The column `delta' denotes the sum of the delta invariants of the singularity.
\item[$\bullet$] The column `deg' denotes a list $(\deg{\textrm{C}}_i)_i$ of the degrees of the the components of ${\textrm{C}}$.
\item[$\bullet$] The column `T5' denotes the number of components of $F$ (either exact or an interval).
\item[$\bullet$] The column `T5R' denotes the number of rational components of $F$ (either exact or an interval).
\item[$\bullet$] We note that not for all root subsystems there exists a cone curve with a corresponding singularity type. These entries are filled with $-$.
\end{itemize}
\newpage
\begin{center}
{\tiny \begin{longtable}
{| l|l|l|l|l|l|l| } \hline index & type & L & delta & deg & T5 & T5R \\
 \hline \hline $ 100 $ & $  A0             $ & $ []                                               $ & $ 0   $ & $ (6)        $ & $ [1, 12]  $ & $ 0        $ \\
 \hline $ 101 $ & $  A1             $ & $ [78]                                             $ & $ 1   $ & $ (6)        $ & $ [1, 10]  $ & $ 0        $ \\
 \hline $ 102 $ & $ 2A1             $ & $ [56, 78]                                         $ & $ 2   $ & $ (6)        $ & $ [1, 8]   $ & $ [0, 1]   $ \\
 \hline $ 103 $ & $  A2             $ & $ [67, 78]                                         $ & $ 1   $ & $ (6)        $ & $ [1, 9]   $ & $ 0        $ \\
 \hline $ 104 $ & $ 3A1             $ & $ [34, 56, 78]                                     $ & $ 3   $ & $ (6)        $ & $ [1, 6]   $ & $ [0, 6]   $ \\
 \hline $ 105 $ & $  A2+ A1         $ & $ [45, 67, 78]                                     $ & $ 2   $ & $ (6)        $ & $ [1, 7]   $ & $ [0, 1]   $ \\
 \hline $ 106 $ & $  A3             $ & $ [56, 67, 78]                                     $ & $ 2   $ & $ (6)        $ & $ [1, 8]   $ & $ [0, 1]   $ \\
 \hline $ 107 $ & $ 4A1             $ & $ [12, 34, 56, 78]                                 $ & $ 4   $ & $ (6)        $ & $ [1, 4]   $ & $ [0, 4]   $ \\
 \hline $ 108 $ & $ 4A1             $ & $ [1145, 1123, 23, 45]                             $ & $ 4   $ & $ (4, 2)     $ & $ 4        $ & $ [3, 4]   $ \\
 \hline $ 109 $ & $  A2+2A1         $ & $ [23, 45, 67, 78]                                 $ & $ 3   $ & $ (6)        $ & $ [1, 5]   $ & $ [0, 5]   $ \\
 \hline $ 110 $ & $ 2A2             $ & $ [34, 45, 67, 78]                                 $ & $ 2   $ & $ (6)        $ & $ [1, 6]   $ & $ [0, 1]   $ \\
 \hline $ 111 $ & $  A3+ A1         $ & $ [34, 56, 67, 78]                                 $ & $ 3   $ & $ (6)        $ & $ [1, 6]   $ & $ [0, 6]   $ \\
 \hline $ 112 $ & $  A4             $ & $ [45, 56, 67, 78]                                 $ & $ 2   $ & $ (6)        $ & $ [1, 7]   $ & $ [0, 1]   $ \\
 \hline $ 113 $ & $  D4             $ & $ [1123, 23, 34, 45]                               $ & $ 3   $ & $ (6)        $ & $ [1, 6]   $ & $ [0, 6]   $ \\
 \hline $ 114 $ & $ 5A1             $ & $ [1145, 1123, 23, 45, 78]                         $ & $ 5   $ & $ (4, 2)     $ & $ 4        $ & $ [3, 4]   $ \\
 \hline $ 115 $ & $  A2+3A1         $ & $ [1123, 23, 45, 67, 78]                           $ & $ 4   $ & $ (6)        $ & $ [1, 3]   $ & $ [0, 3]   $ \\
 \hline $ 116 $ & $ 2A2+ A1         $ & $ [12, 34, 45, 67, 78]                             $ & $ 3   $ & $ (6)        $ & $ [1, 4]   $ & $ [0, 4]   $ \\
 \hline $ 117 $ & $  A3+2A1         $ & $ [12, 34, 56, 67, 78]                             $ & $ 4   $ & $ (6)        $ & $ [1, 4]   $ & $ [0, 4]   $ \\
 \hline $ 118 $ & $  A3+2A1         $ & $ [1145, 1123, 23, 45, 56]                         $ & $ 4   $ & $ (4, 2)     $ & $ 4        $ & $ [3, 4]   $ \\
 \hline $ 119 $ & $  A3+ A2         $ & $ [23, 34, 56, 67, 78]                             $ & $ 3   $ & $ (6)        $ & $ [1, 5]   $ & $ [0, 5]   $ \\
 \hline $ 120 $ & $  A4+ A1         $ & $ [23, 45, 56, 67, 78]                             $ & $ 3   $ & $ (6)        $ & $ [1, 5]   $ & $ [0, 5]   $ \\
 \hline $ 121 $ & $  A5             $ & $ [34, 45, 56, 67, 78]                             $ & $ 3   $ & $ (6)        $ & $ [1, 6]   $ & $ [0, 6]   $ \\
 \hline $ 122 $ & $  D4+ A1         $ & $ [1123, 23, 34, 45, 78]                           $ & $ 4   $ & $ (6)        $ & $ [1, 4]   $ & $ [0, 4]   $ \\
 \hline $ 123 $ & $  D5             $ & $ [1123, 23, 34, 45, 56]                           $ & $ 3   $ & $ (6)        $ & $ [1, 5]   $ & $ [0, 5]   $ \\
 \hline $ 124 $ & $ 6A1             $ & $ [1567, 1347, 1127, 12, 34, 56]                   $ & $ 6   $ & $ (2, 2, 2)  $ & $ 3        $ & $ 3        $ \\
 \hline $ 125 $ & $  A2+4A1         $ & $ [1145, 1123, 23, 45, 67, 78]                     $ & $ 5   $ & $ (4, 2)     $ & $ 4        $ & $ [3, 4]   $ \\
 \hline $ 126 $ & $ 2A2+2A1         $ & $ [1123, 12, 23, 45, 67, 78]                       $ & $ 4   $ & $ (6)        $ & $ [1, 2]   $ & $ [0, 2]   $ \\
 \hline $ 127 $ & $ 3A2             $ & $ [1456, 1123, 12, 23, 45, 56]                     $ & $ 3   $ & $ (6)        $ & $ [1, 3]   $ & $ [0, 3]   $ \\
 \hline $ 128 $ & $  A3+3A1         $ & $ [1145, 1123, 23, 45, 56, 78]                     $ & $ 5   $ & $ (4, 2)     $ & $ 4        $ & $ [3, 4]   $ \\
 \hline $ 129 $ & $  A3+ A2+ A1     $ & $ [1123, 12, 34, 56, 67, 78]                       $ & $ 4   $ & $ (6)        $ & $ [1, 3]   $ & $ [0, 3]   $ \\
 \hline $ 130 $ & $ 2A3             $ & $ [12, 23, 34, 56, 67, 78]                         $ & $ 4   $ & $ (6)        $ & $ [1, 4]   $ & $ [0, 4]   $ \\
 \hline $ 131 $ & $ 2A3             $ & $ [1567, 1145, 1127, 12, 56, 78]                   $ & $ 4   $ & $ (4, 2)     $ & $ 4        $ & $ [3, 4]   $ \\
 \hline $ 132 $ & $  A4+2A1         $ & $ [1123, 23, 45, 56, 67, 78]                       $ & $ 4   $ & $ (6)        $ & $ [1, 3]   $ & $ [0, 3]   $ \\
 \hline $ 133 $ & $  A4+ A2         $ & $ [12, 23, 45, 56, 67, 78]                         $ & $ 3   $ & $ (6)        $ & $ [1, 4]   $ & $ [0, 4]   $ \\
 \hline $ 134 $ & $  A5+ A1         $ & $ [12, 34, 45, 56, 67, 78]                         $ & $ 4   $ & $ (6)        $ & $ [1, 4]   $ & $ [0, 4]   $ \\
 \hline $ 135 $ & $  A5+ A1         $ & $ [1145, 1123, 12, 23, 45, 56]                     $ & $ 4   $ & $ (4, 2)     $ & $ 4        $ & $ [3, 4]   $ \\
 \hline $ 136 $ & $  A6             $ & $ [23, 34, 45, 56, 67, 78]                         $ & $ 3   $ & $ (6)        $ & $ [1, 5]   $ & $ [0, 5]   $ \\
 \hline $ 137 $ & $  D4+2A1         $ & $ [1145, 1123, 23, 45, 56, 67]                     $ & $ 5   $ & $ (4, 2)     $ & $ 4        $ & $ [3, 4]   $ \\
 \hline $ 138 $ & $  D4+ A2         $ & $ [1123, 23, 34, 45, 67, 78]                       $ & $ 4   $ & $ (6)        $ & $ [1, 3]   $ & $ [0, 3]   $ \\
 \hline $ 139 $ & $  D5+ A1         $ & $ [1123, 23, 34, 45, 56, 78]                       $ & $ 4   $ & $ (6)        $ & $ [1, 3]   $ & $ [0, 3]   $ \\
 \hline $ 140 $ & $  D6             $ & $ [1123, 23, 34, 45, 56, 67]                       $ & $ 4   $ & $ (6)        $ & $ [1, 4]   $ & $ [0, 4]   $ \\
 \hline $ 141 $ & $  E6             $ & $ [1123, 12, 23, 34, 45, 56]                       $ & $ 3   $ & $ (6)        $ & $ [1, 4]   $ & $ [0, 4]   $ \\
 \hline $ 142 $ & $ 7A1             $ & $ [278, 1567, 1347, 1127, 12, 34, 56]              $ & $ 7   $ & $ -          $ & $ -        $ & $ -        $ \\
 \hline $ 143 $ & $ 3A2+ A1         $ & $ [1456, 1123, 12, 23, 45, 56, 78]                 $ & $ 4   $ & $ (6)        $ & $ 1        $ & $ [0, 1]   $ \\
 \hline $ 144 $ & $  A3+4A1         $ & $ [278, 1347, 1127, 12, 34, 56, 78]                $ & $ 6   $ & $ (2, 2, 2)  $ & $ 3        $ & $ 3        $ \\
 \hline $ 145 $ & $  A3+ A2+2A1     $ & $ [1145, 1123, 12, 23, 45, 67, 78]                 $ & $ 5   $ & $ (4, 2)     $ & $ 4        $ & $ [3, 4]   $ \\
 \hline $ 146 $ & $ 2A3+ A1         $ & $ [1678, 1456, 1347, 1123, 12, 45, 78]             $ & $ 5   $ & $ (4, 2)     $ & $ 4        $ & $ [3, 4]   $ \\
 \hline $ 147 $ & $  A4+ A2+ A1     $ & $ [1123, 12, 23, 45, 56, 67, 78]                   $ & $ 4   $ & $ (6)        $ & $ [1, 2]   $ & $ [0, 2]   $ \\
 \hline $ 148 $ & $  A4+ A3         $ & $ [1123, 12, 23, 34, 56, 67, 78]                   $ & $ 4   $ & $ (6)        $ & $ [1, 3]   $ & $ [0, 3]   $ \\
 \hline $ 149 $ & $  A5+2A1         $ & $ [1145, 1123, 12, 23, 45, 56, 78]                 $ & $ 5   $ & $ (4, 2)     $ & $ 4        $ & $ [3, 4]   $ \\
 \hline $ 150 $ & $  A5+ A2         $ & $ [1456, 1123, 12, 23, 45, 56, 67]                 $ & $ 4   $ & $ (6)        $ & $ [1, 3]   $ & $ [0, 3]   $ \\
 \hline $ 151 $ & $  A6+ A1         $ & $ [1123, 12, 34, 45, 56, 67, 78]                   $ & $ 4   $ & $ (6)        $ & $ [1, 3]   $ & $ [0, 3]   $ \\
 \hline $ 152 $ & $  A7             $ & $ [12, 23, 34, 45, 56, 67, 78]                     $ & $ 4   $ & $ (6)        $ & $ [1, 4]   $ & $ [0, 4]   $ \\
 \hline $ 153 $ & $  A7             $ & $ [1567, 1145, 1127, 12, 23, 56, 78]               $ & $ 4   $ & $ (4, 2)     $ & $ 4        $ & $ [3, 4]   $ \\
 \hline $ 154 $ & $  D4+3A1         $ & $ [1567, 1347, 1127, 12, 34, 56, 78]               $ & $ 6   $ & $ (2, 2, 2)  $ & $ 3        $ & $ 3        $ \\
 \hline $ 155 $ & $  D4+ A3         $ & $ [1567, 1347, 1145, 1127, 12, 56, 78]             $ & $ 5   $ & $ (4, 2)     $ & $ 4        $ & $ [3, 4]   $ \\
 \hline $ 156 $ & $  D5+2A1         $ & $ [1145, 1123, 23, 45, 56, 67, 78]                 $ & $ 5   $ & $ (4, 2)     $ & $ 4        $ & $ [3, 4]   $ \\
 \hline $ 157 $ & $  D5+ A2         $ & $ [1123, 12, 23, 34, 45, 67, 78]                   $ & $ 4   $ & $ (6)        $ & $ [1, 2]   $ & $ [0, 2]   $ \\
 \hline $ 158 $ & $  D6+ A1         $ & $ [1145, 1123, 12, 23, 45, 56, 67]                 $ & $ 5   $ & $ (4, 2)     $ & $ 4        $ & $ [3, 4]   $ \\
 \hline $ 159 $ & $  D7             $ & $ [1123, 23, 34, 45, 56, 67, 78]                   $ & $ 4   $ & $ (6)        $ & $ [1, 3]   $ & $ [0, 3]   $ \\
 \hline $ 160 $ & $  E6+ A1         $ & $ [1123, 12, 23, 34, 45, 56, 78]                   $ & $ 4   $ & $ (6)        $ & $ [1, 2]   $ & $ [0, 2]   $ \\
 \hline $ 161 $ & $  E7             $ & $ [1123, 12, 23, 34, 45, 56, 67]                   $ & $ 4   $ & $ (6)        $ & $ [1, 4]   $ & $ [0, 4]   $ \\
 \hline $ 162 $ & $ 8A1             $ & $ [308, 278, 1567, 1347, 1127, 12, 34, 56]         $ & $ 8   $ & $ -          $ & $ -        $ & $ -        $ \\
 \hline $ 163 $ & $ 4A2             $ & $ [1123, 1345, 1156, 1258, 1367, 1247, 1468, 1178] $ & $ 4   $ & $ (6)        $ & $ 0        $ & $ [0, 0]   $ \\
 \hline $ 164 $ & $ 2A3+2A1         $ & $ [308, 278, 1567, 12, 23, 34, 56, 67]             $ & $ 6   $ & $ (2, 2, 2)  $ & $ 3        $ & $ 3        $ \\
 \hline $ 165 $ & $ 2A4             $ & $ [278, 1678, 12, 23, 34, 45, 67, 78]              $ & $ 4   $ & $ (6)        $ & $ [1, 2]   $ & $ [0, 2]   $ \\
 \hline $ 166 $ & $  A5+ A2+ A1     $ & $ [1678, 1145, 1123, 12, 23, 45, 67, 78]           $ & $ 5   $ & $ (4, 2)     $ & $ 4        $ & $ [3, 4]   $ \\
 \hline $ 167 $ & $  A7+ A1         $ & $ [234, 1145, 1123, 12, 23, 56, 67, 78]            $ & $ 5   $ & $ (4, 2)     $ & $ 4        $ & $ [3, 4]   $ \\
 \hline $ 168 $ & $  A8             $ & $ [1567, 1123, 12, 23, 34, 56, 67, 78]             $ & $ 4   $ & $ (6)        $ & $ [1, 3]   $ & $ [0, 3]   $ \\
 \hline $ 169 $ & $  D4+4A1         $ & $ [278, 1567, 1347, 1127, 12, 34, 56, 78]          $ & $ 7   $ & $ -          $ & $ -        $ & $ -        $ \\
 \hline $ 170 $ & $ 2D4             $ & $ [234, 278, 12, 23, 34, 56, 67, 78]               $ & $ 6   $ & $ (2, 2, 2)  $ & $ 3        $ & $ 3        $ \\
 \hline $ 171 $ & $  D5+ A3         $ & $ [278, 1567, 12, 23, 34, 56, 67, 78]              $ & $ 5   $ & $ (4, 2)     $ & $ 4        $ & $ [3, 4]   $ \\
 \hline $ 172 $ & $  D6+2A1         $ & $ [278, 1347, 1127, 12, 34, 45, 56, 78]            $ & $ 6   $ & $ (2, 2, 2)  $ & $ 3        $ & $ 3        $ \\
 \hline $ 173 $ & $  D8             $ & $ [278, 12, 23, 34, 45, 56, 67, 78]                $ & $ 5   $ & $ (4, 2)     $ & $ 4        $ & $ [3, 4]   $ \\
 \hline $ 174 $ & $  E6+ A2         $ & $ [1456, 1123, 12, 23, 45, 56, 67, 78]             $ & $ 4   $ & $ (6)        $ & $ 1        $ & $ [0, 1]   $ \\
 \hline $ 175 $ & $  E7+ A1         $ & $ [1145, 1123, 12, 23, 45, 56, 67, 78]             $ & $ 5   $ & $ (4, 2)     $ & $ 4        $ & $ [3, 4]   $ \\
 \hline $ 176 $ & $  E8             $ & $ [1123, 12, 23, 34, 45, 56, 67, 78]               $ & $ 4   $ & $ (6)        $ & $ 1        $ & $ [0, 1]   $ \\
 \hline
\end{longtable}
 }
\end{center}
\end{tab}

\section{Algorithm: radical decomposition of equidimensional ideals}
The following algorithm turned out to be much faster than other radical decomposition algorithms which were available to us. We compute the radical decomposition of equidimensional ideals and is used by the next two algorithms. The idea of this algorithm is to project an algebraic subset and compute the radical components of the projection. For each radical component we intersect the cone over this component with the original algebraic subset. Then we do a radical decomposition on this intersection.
\begin{algorithm}\label{alg:f2_equidim_radical_decomp}\textbf{\textrm{(radical decomposition of equidimensional ideals)}}{\small
\begin{itemize}
\item[$\bullet$]\textbf{function:~} EquidimRadicalDecomposition$(~I,~x_k~)$
\item[]\textit{input:~} An ideal $I$ in some polynomial ring $R={\textbf{F}}[x_0,\ldots,x_n]$ such that the variety of each primary component has equal dimension. A generator $x_k \in R$.
\item[]\textit{output:~} The radical primary decomposition of $I$. First the algorithm projects $V(I)$ from $x_k$. Then computes the radical components of the projection. We return the intersection with $V(I)$ of the cones over the projected radical components. If the projection was not one-to-one then the radical decomposition may contain artifacts.
\item[] $I_e:=$Eliminate$(~I,~[x_k]~)$
\item[] $(I_{ei})_i:=$RadicalDecomposition$(~I_e~)$
\item[] $(I_i)_i:=(~$RadicalDecomposition$(~I+I_{ei}~)~)_i$
\item[]\textbf{return}~ $(I_i)_i$
\item[]\textbf{end function}
\end{itemize} }
\end{algorithm}

\begin{proposition}\label{prop:}\textbf{\textrm{(radical decomposition of equidimensional ideals)}}

\textbf{a)} We have that Algorithm~\ref{alg:f2_equidim_radical_decomp} is correct.

\begin{proof}
 See \cite{nls2}, chapter 8, section 6, proposition 113, page 164.
 \end{proof}
\end{proposition}

\section{Algorithm: analyze tangent developable surface of dual of cone curve}
This algorithm computes the components of the singular locus of the tangent developable of the dual cone curve. From Proposition~\ref{prop:f2_cone_curve_dual_surface} it follows that these singular components define families of bitangent planes.
\newpage
\begin{algorithm}\label{alg:f2_dual_tangent_developable}\textbf{\textrm{(tangent developable of dual of cone curve)}}{\small
\begin{itemize}
\item[$\bullet$]\textbf{function:~} AnalyseDualSurface$(~W~)$
\item[]\textit{input:~} A homogeneous polynomial $W(x,y,z) \in {\textbf{F}}[x,y,z]$ of degree $6$ with weights $(2:1:1)$, which is the equation of the weighted cone curve of a cone curve ${\textrm{C}}$.
\item[]\textit{output:~} The radical ideals of the irreducible components of the singular locus of the tangent surface of the second associated curve of ${\textrm{C}}$ (see Proposition~\ref{prop:f2_cone_curve_dual_surface}).
\item[] $H := k x + l y^2 + m y z + n z^2$
\item[] $I:={\textbf{F}}\langle ~ \partial_xH \partial_yW - \partial_yH \partial_xW,~ \partial_xH \partial_zW - \partial_zH \partial_xW,~ \partial_yH \partial_zW - \partial_zH \partial_yW ~\rangle $
\item[] $R:={\textbf{Q}}\langle ~y,~z~\rangle $
\item[] $I_s:=(I:R)$
\item[] $I_e:=$Eliminate$(~I_s,~[x,y,z]~)$
\item[] $g:=$NonLinearIrreducibleFactor(~GeneratorOfPrincipalIdeal(~$I_e$~)~)
\item[] $S:={\textbf{F}}\langle ~\partial_kg,~\partial_lg,~\partial_mg,~\partial_ng~\rangle $
\item[]\textbf{return}~ EquidimRadicalDecomposition$(~S,~k~)$ (Algorithm~\ref{alg:f2_equidim_radical_decomp})
\item[]\textbf{end function}
\end{itemize} }
\end{algorithm}

\begin{proposition}\label{prop:f2_dual_tangent_developable}\textbf{\textrm{(properties of analyze dual surface algorithm)}}

\textbf{a)} We have that Algorithm~\ref{alg:f2_dual_tangent_developable} is correct.

\begin{proof}
 See \cite{nls2}, chapter 8, section 7, proposition 114, page 165.
 \end{proof}
\end{proposition}
The workings and output of the above algorithm is presented for a concrete elliptic cone curve with three cusps. There are three irreducible components of the singular locus of the tangent developable surface of the dual cone curve, which correspond to the three components of the reduced bitangent correspondence.

 In this example the three components in the singular locus are dual to the families of tangent planes of three cubic cones. These are cones in the linear series of cubics with the cone curve as base locus.

 From Theorem~\ref{thm:f2_coc_cls} we know that the upper bound for the number of components of the reduced bitangent correspondence is three.
\begin{example}\label{ex:f2_dual_tangent_developable}\textrm{\textbf{(analyze dual surface algorithm)}}
 Let ${\textrm{C}} \subset {\textbf{P}}^3(s:t:u:v)$ be a cone curve defined by the complete intersection of a quadric cone ${\textrm{Q}}$ and cubic surface ${\textrm{U}}$. We use the notation as in Algorithm~\ref{alg:f2_dual_tangent_developable}. Let ${\textrm{W}}: W(x,y,z)=x^3 + y^4 z^2 - 2 y^3 z^3 + y^2 z^4=0 \subset{\textbf{P}}(2:1:1) $ be the weighted cone curve of ${\textrm{C}}$.

 We find that ${\textrm{W}}$ has three $A_2$ singularities at $(0:0:1)$, $(0:1:0)$ and $(0:1:1)$.

 We call AnalyseDualSurface$(W)$.

 We find that $I_e={\textbf{Q}}\langle ~ n^3 l^3 (l + m + n)^3 g(k,l,m,n) ~\rangle $ for some polynomial $g(k,l,m,n)$ of degree $9$. Note that the powers of the linear factors are related to the singularities.

 We have that ${\textrm{T}}: g(k,l,m,n)=0 \subset{\textbf{P}}^{3*}$ is the tangent developable of the second associated curve of ${\textrm{C}}$. The singular locus on ${\textrm{T}}$ is defined by the zeroset of $S:={\textbf{Q}}\langle ~\partial_kg,~\partial_lg,~\partial_mg,~\partial_ng~\rangle )$.

 From Proposition~\ref{prop:f2_cone_curve_dual_surface}.e) it follows that $\deg{\textrm{T}}=2g+10-\beta-2\gamma=9$ where $\beta=3,\gamma=0$ and $g=p_g({\textrm{C}})=1$.

 We obtain the following irreducible radical components $S_1,\ldots, S_5$ (the multiplicities are omitted) of $S$:
\begin{itemize}\addtolength{\itemsep}{1pt}
\item[$\bullet$] $ S_1 = {\textbf{Q}} \langle ~ l^2 m^2 + l m^3 + \frac{1}{16} m^4 + 32 l^3 n + 50 l^2 m n + \frac{41}{2} l m^2 n + m^3 n + 61 l^2 n^2 + 50 l m n^2 + m^2 n^2 + 32 l n^3 ,~ k^3 + \frac{3}{4} l m^2 + \frac{3}{8} m^3 + 24 l^2 n + \frac{51}{2} l m n + \frac{3}{4} m^2 n + 24 l n^2 ~\rangle $,
\item[$\bullet$] $S_2 = {\textbf{Q}}\langle ~ m^2 - 4 l n ,~ k ~\rangle $,
\item[$\bullet$] $S_3 = {\textbf{Q}}\langle ~ k^3 - \frac{27}{4} m^2 n ,~ l + m ~\rangle $,
\item[$\bullet$] $S_4 = {\textbf{Q}}\langle ~ k^3 - \frac{27}{4} m^2 n - 27 m n^2 - 27 n^3 ,~ l - n ~\rangle $, and
\item[$\bullet$] $S_5 = {\textbf{Q}}\langle ~ k^3 - \frac{27}{4} l n^2 ,~ m + n ~\rangle $.
\end{itemize}

 The component $V(S_1)$ is the cuspidal component dual to the 3-osculating planes of ${\textrm{C}}$.

 The component $V(S_2)$ is a plane conic which is dual to the tritangent family on ${\textrm{C}}$ (thus the planes tangent to the quadric cone ${\textrm{Q}}$). These planes are dual to a conic.

 The components $V(S_3)$, $V(S_4)$ and $V(S_5)$ define irreducible components of the reduced bitangent correspondence. Since the components are planar it follows that the dual planes all go through a point. It follows that these three families of bitangent planes are defined by the tangent planes of cubic cones in the linear series $\Lambda={\textbf{Q}}\langle {\textrm{U}},s{\textrm{Q}},t{\textrm{Q}},u{\textrm{Q}},v{\textrm{Q}}\rangle $. It can be shown that there are at most three cubic cones in $\Lambda$.

 From Proposition~\ref{prop:f2_cone_curve_bound} it follows that the upper bound for the number of bitangent families is $n=3$, which in this example is reached.

 From Proposition~\ref{prop:f2_cone_curve_dual_surface}.c) it follows that indeed $\deg S_3\underset{}{\sqcup} S_4\underset{}{\sqcup} S_5 = 9$ where $\beta=3$ and $g=p_g({\textrm{C}})=1$.

 We have that $V(S_1)=1$ and $p_gV(S_i)=1$ for $i>2$. From the proof of Proposition~\ref{prop:f2_cone_curve_correspondence_genus} it follows that there is an unramified $2:1$ morphism defined by ${\textrm{T}}$ from $V(S_1)$ onto $V(S_i)$.

\end{example}

\section{Algorithm: analyze tangent developable surface of cone curve}
This algorithm computes the components of the singular locus of the tangent developable the cone curve itself. Recall from Proposition~\ref{prop:f2_cone_curve_dual_surface} that there is a bijection between the components of the singular loci of the tangent developables of the cone curve and dual cone curve.
\begin{algorithm}\label{alg:f2_analyze_developable_surface}\textbf{\textrm{(analyze developable surface of cone curve)}}{\small
\begin{itemize}
\item[$\bullet$]\textbf{function:~} AnalyseDevelopableSurface$(~U~)$
\item[]\textit{input:~} A homogeneous polynomial $U(s,t,u,v) \in {\textbf{C}}[s,t,u,v]$ of degree $3$.
\item[]\textit{output:~} The one dimensional irreducible components of the singular locus of the tangent developable surface of the cone curve $V(Q,U)\subset {\textbf{P}}^3$ where $Q=u^2 - tv$.
\item[] $Q:=u^2 - tv$
\item[] $p:=(a:b:c:d)$
\item[] $TQ:=s \partial_s Q + t \partial_t Q + u \partial_u Q + v \partial_v Q$
\item[] $TU:=s \partial_s U + t \partial_t U + u \partial_u U + v \partial_v U$
\item[] $I:={\textbf{F}}\langle ~TQ(p),~TU(p),~Q(p),~U(p)~\rangle $
\item[] $R:={\textbf{F}}\langle ~a,~b,~c,~d~\rangle $
\item[] $I_s:=(I:R)$
\item[] $I_e:=$Eliminate$(~I_s,~[~a,~b,~c,~d~]~)$
\item[] $g:=$NonLinearIrreducibleFactor(~GeneratorOfPrincipalIdeal(~$I_e$~)~)
\item[] $M:={\textbf{F}}\langle ~\partial_sg,~\partial_tg,~\partial_ug,~\partial_vg~\rangle $
\item[]\textbf{return}~ EquidimRadicalDecomposition$(M)$ (Algorithm~\ref{alg:f2_equidim_radical_decomp})
\item[]\textbf{end function}
\end{itemize} }
\end{algorithm}

\begin{proposition}\label{prop:}\textbf{\textrm{(properties of AnalyseDevelopableSurface algorithm)}}

\textbf{a)} We have that Algorithm~\ref{alg:f2_analyze_developable_surface} is correct.

\begin{proof}
 See \cite{nls2}, chapter 8, section 9, proposition 118, page 170.
 \end{proof}
\end{proposition}
We analyze the singular locus of the tangent developable of the same elliptic cone curve which we analyzed in the previous example. We go step by step through the algorithm and show the output.
\begin{example}\label{ex:f2_analyze_developable_surface}\textrm{\textbf{(analyze developable surface algorithm)}}
  We us the same notation as in Algorithm~\ref{alg:f2_analyze_developable_surface}. Let ${\textrm{U}}: U(s,t,u,v)=s^3+u^2t-2u^3+u^2v \subset {\textbf{P}}^3$ be a cubic surface.  Let ${\textrm{Q}}: u^2-tv=0$ be a quadric cone.  Let ${\textrm{C}} \subset {\textbf{P}}^3$ be the cone curve defined as the complete intersection of ${\textrm{U}}$ and ${\textrm{Q}}$.

 Note that the weighted cone curve (see Definition~\ref{def:f2_weighted_cone_curve}) of ${\textrm{C}}$ is ${\textrm{W}}: W(x,y,z)=x^3 + y^4 z^2 - 2 y^3 z^3 + y^2 z^4=0 \subset{\textbf{P}}(2:1:1) $ as in Example~\ref{ex:f2_dual_tangent_developable}.

 We have that ${\textrm{C}}$ has $3$ $A_2$ singularities and $p_g({\textrm{C}})=1$.

 We call AnalyseDevelopableSurface$(U)$ at Algorithm~\ref{alg:f2_analyze_developable_surface}.

 We find that $I_e={\textbf{Q}}\langle  v^2  t^2 (t - 2u + v)^2 g(s,t,u,v)\rangle $ for some polynomial $g(s,t,u,v)$ of degree $9$. Note that each singularity results in a power of a linear factor.

 We have that ${\textrm{Z}}: g(s,t,u,v)=0 \subset{\textbf{P}}^{3}$ is the tangent developable of ${\textrm{C}}$. The singular locus on ${\textrm{Z}}$ is defined by the zeroset of $M:={\textbf{Q}}\langle ~\partial_sg,~\partial_tg,~\partial_ug,~\partial_vg~\rangle )$.

 From Proposition~\ref{prop:f2_cone_curve_dual_surface}.e) it follows that indeed $\deg{\textrm{Z}}=2g+10-\beta-2\gamma=9$ where $\beta=3,\gamma=0$ and $g=p_g({\textrm{C}})=1$.

 We obtain the following irreducible radical components $M_1,\ldots, M_5$ (the multiplicities are omitted) of $M$:
\begin{itemize}\addtolength{\itemsep}{1pt}
\item[$\bullet$] $ M_1 = {\textbf{Q}}\langle ~ s^3 + t^2 v - 2 t u v + t v^2 ,~ u^2 - t v ~\rangle $,
\item[$\bullet$] $ M_2 = {\textbf{Q}}\langle ~ t u^2 - u^3 + \frac{1}{8} t^2 v - \frac{5}{4} t u v + u^2 v + \frac{1}{8} t v^2 ,~ s ~\rangle  $,
\item[$\bullet$] $ M_3 = {\textbf{Q}}\langle ~ u^4 + \frac{27}{16} s^3 v - \frac{3}{2} u^3 v + \frac{3}{4} u^2 v^2 - \frac{1}{8} u v^3 ,~ t - 2 u ~\rangle $,
\item[$\bullet$] $ M_4 = {\textbf{Q}}\langle ~ u^4 + \frac{27}{16} s^3 v - \frac{5}{2} u^3 v + \frac{9}{4} u^2 v^2 - \frac{7}{8} u v^3 + \frac{1}{8} v^4 ,~ t - v ~\rangle  $, and
\item[$\bullet$] $ M_5 = {\textbf{Q}}\langle ~ s^3 t - \frac{1}{27} t^3 v + \frac{1}{9} t^2 v^2 - \frac{1}{9} t v^3 + \frac{1}{27} v^4 ,~ u - \frac{1}{2} v ~\rangle  $.
\end{itemize}

 We have that $V(M_1)$ is the cuspidal curve defined by ${\textrm{C}}$ on the tangent developable ${\textrm{Z}}$.

 We have that $V(M_2)$ are the intersections of the tangent lines at three points lying on a line through the vertex. The plane containing these tangent lines is a tritangent plane, tangent to ${\textrm{Q}}$. We have that $V(M_2)$ has genus zero and $V(M_2)$ defines the tritangent correspondence.

 We have that $V(M_3), V(M_4)$ and $V(M_5)$ define the irreducible components of the reduced bitangent correspondence. A generic point on $\underset{i>2}{\sqcup}V(M_i)$ is the intersection of two tangent lines of ${\textrm{C}}$. The plane containing these two lines is a bitangent plane of ${\textrm{C}}$. We have that the genus $p_g(M_i)=0$ for $i=3,4,5$.

 In the images below we see an affine chart of the cone curve ${\textrm{C}}$ lying on a quadric cone ${\textrm{Q}}$ and its tangent developable ${\textrm{Z}}$. The singular components $V(M_3), V(M_4)$ and $V(M_5)$ are colored white, red and blue  respectively . The three images are at different scalings. The right most image is from a different viewpoint then the first two.
\begin{center}
 \begin{tabular}
{ccc} {\includegraphics[width=4cm,height=4cm]{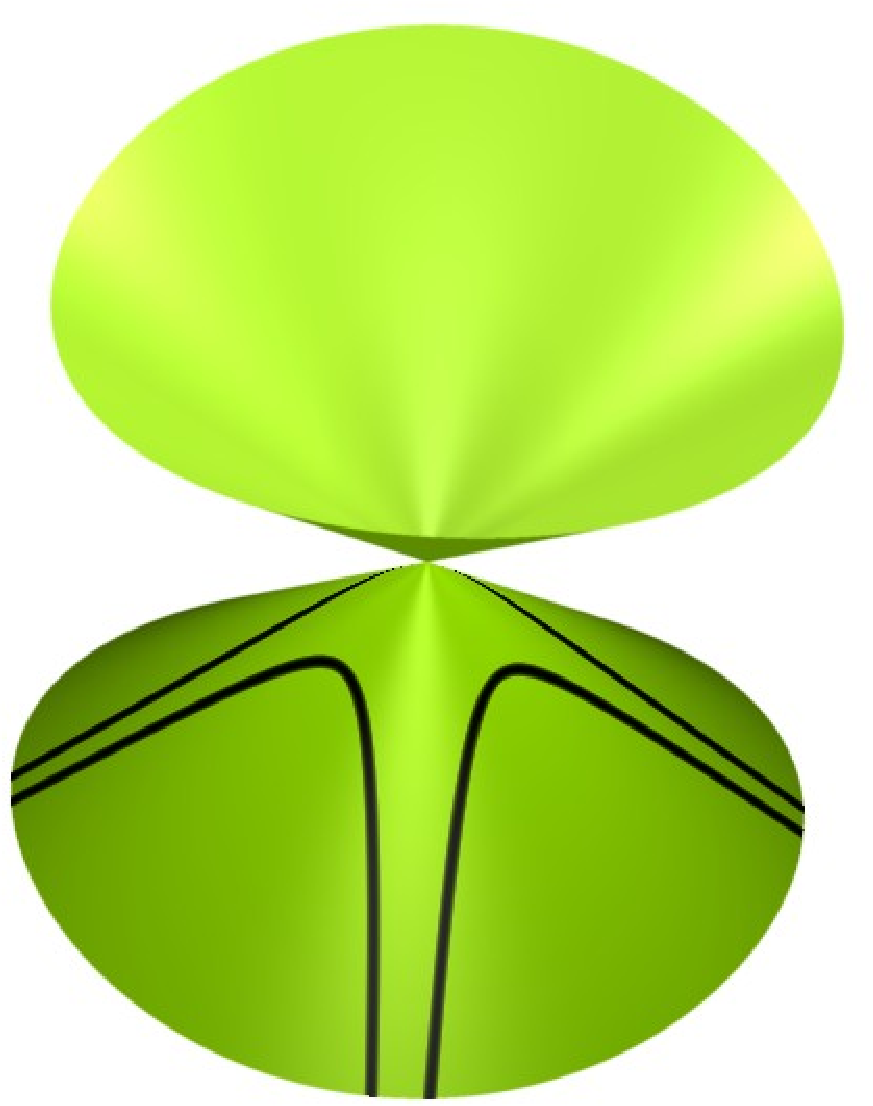}}& {\includegraphics[width=4cm,height=4cm]{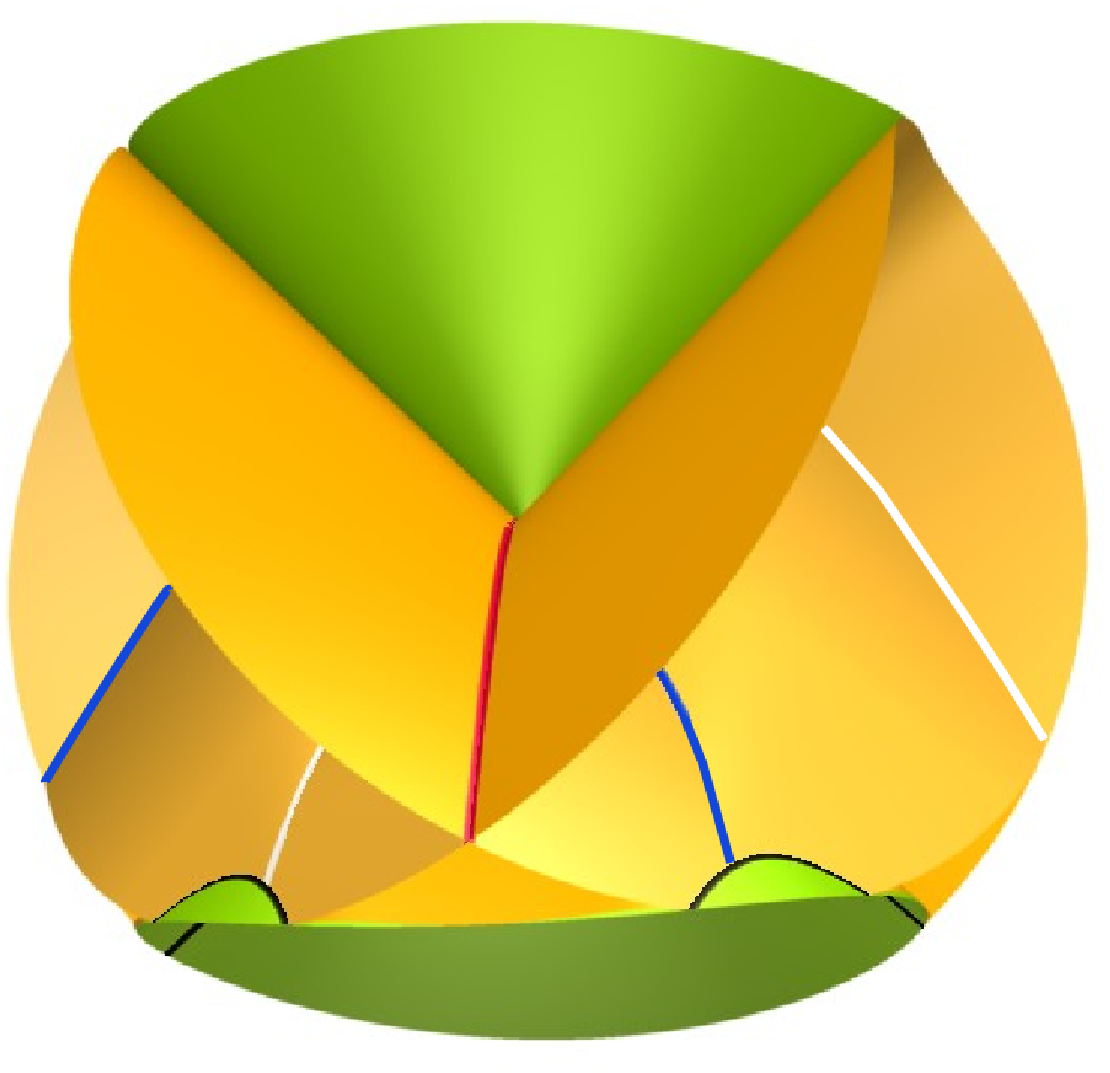}}& {\includegraphics[width=4cm,height=4cm]{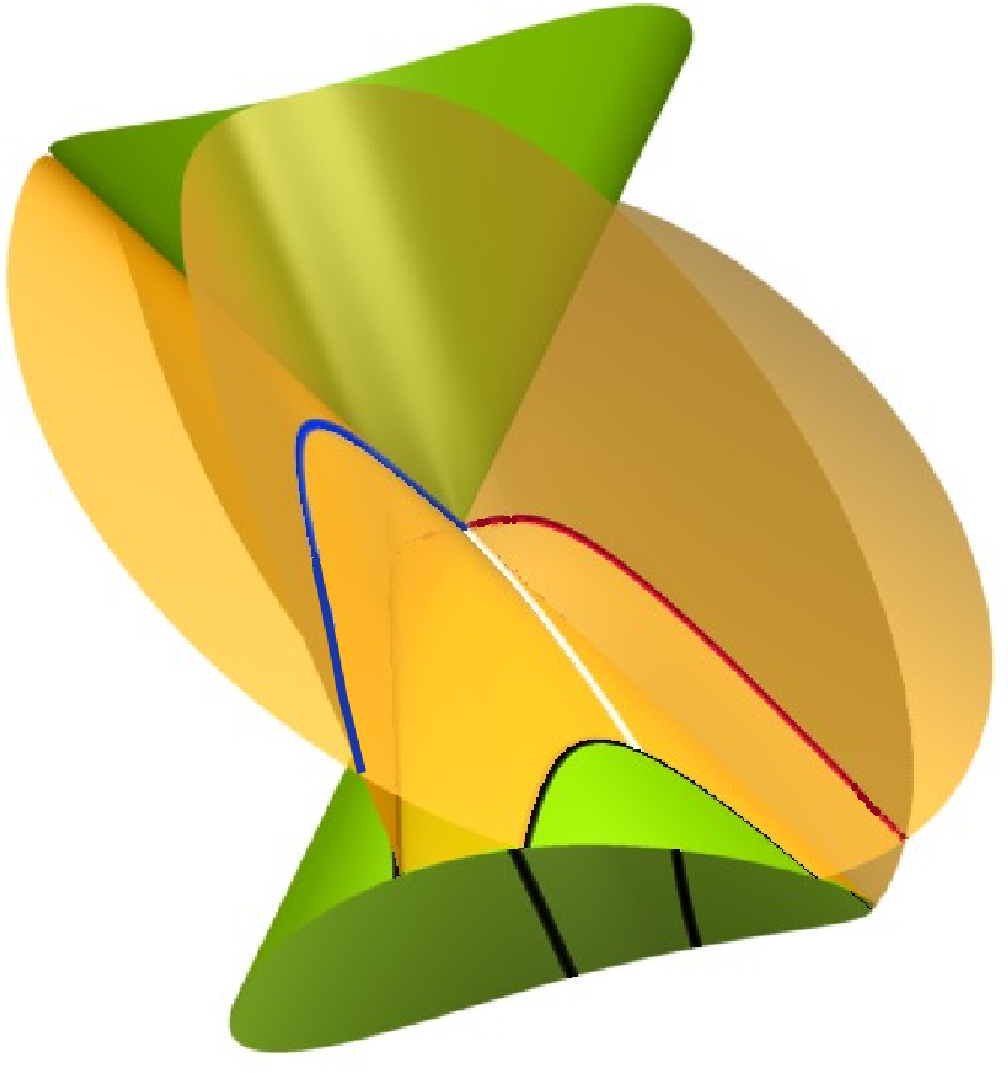}}
\end{tabular}

\end{center} The images were made using \cite{surfex}.

\end{example}

\section{Example non-fibration family on a weak degree one Del Pezzo surface}
Recall that in the motivation section of the introduction we considered an example of a minimal non-fibration family of the plane. This family was defined by the tangent lines of the unit circle.

 In the following example we consider a minimal non-fibration family of a weak degree one Del Pezzo surface. Instead of tangent lines of the unit circle we now consider the pull back of bitangent planes of the elliptic cone curve with three cusps as discussed in the previous two sections.
\begin{example}\label{ex:f2_T5}\textrm{\textbf{(non-fibration family on a weak degree one Del Pezzo surface)}}
 Let ${\textrm{W}}: W(x,y,z)= w^2+z^3+x^2y^2(x-y)^2= 0 \subset {\textbf{P}}(1:1:2)$ be as in Example~\ref{ex:f2_dual_tangent_developable}: a weighted cone curve with $3A_2$ singularities. Let ${\textrm{Y}}': w^2 + W(x,y,z) = 0 \subset {\textbf{P}}(1:1:2:3)$ be the weighted homogenous equation of a weak degree one Del Pezzo surface.

 We consider the following embedding in homogeneous space (see Proposition~\ref{prop:f2_DP1_ring}): \[ \ensuremath{i:{\textrm{Y}}'\rightarrow {\textrm{Y}}\subset{\textbf{P}}^6,\quad (x:y:z:w)\mapsto (x^3:x^2y:xy^2:y^3:xz:yz:w)=(a:b:c:d:e:f:g)}. \] Let ${\textrm{X}}$ be the minimal resolution of singularities of ${\textrm{Y}}$ (see \cite{har1}). Let $K$ be the canonical divisor class of ${\textrm{X}}$.

 From Proposition~\ref{prop:f2_cone_curve_dp_function} it follows that ${\textrm{Y}}$ has $3A_2$ singularities and $\ensuremath{{\textrm{Y}}\stackrel{\varphi_{-2K}}{\rightarrow}{\textbf{P}}(2:1:1)}$ has ${\textrm{W}}$ as branch locus. We can represent ${\textrm{Y}}$ as a polarized surface $({\textrm{X}},D)$  such that  $D=-3K$ and $\ensuremath{{\textrm{X}}\stackrel{\varphi_{D}}{\rightarrow}{\textrm{Y}}\subset{\textbf{P}}^6}$. It follows that ${\textrm{Y}}$ has degree $D^2=9$ and the sectional genus (the geometric genus of a generic hyperplane section) equals $p_a(D)=\frac{1}{2}D(D+K)+1=4$. Let ${\textrm{Q}}:t^2-su \subset {\textbf{P}}^3$ be the quadric cone. Let $\ensuremath{\mu:{\textrm{W}}\subset{\textbf{P}}(1:1:2)\rightarrow {\textrm{C}}\subset{\textrm{Q}},\quad (x:y:z)\mapsto (x^2:xy:y^2:z)=(s:t:u:v)}$ be the cone curve isomorphism.

 We obtain a map from ${\textrm{Y}}$ to ${\textrm{Q}}$ wich commutes with $\mu$ and has the cone curve ${\textrm{C}}$ as branch locus: \[ \ensuremath{q:{\textrm{Y}}\rightarrow {\textrm{Q}},\quad (a:b:c:d:e:f:g)\mapsto (a:b:c:e)=(s:t:u:v)}. \]

 We project ${\textrm{Y}}$ to ${\textbf{P}}^3$: \[ \ensuremath{p:{\textrm{Y}}\subset{\textbf{P}}^6\rightarrow {\textrm{Z}}\subset{\textbf{P}}^3,\quad (a:b:c:d:e:f:g)\mapsto (a:d:g:e)=(s:t:u:v)}. \] We compute the image $p({\textrm{Y}})$ using Gr\"obner basis (see \cite{Schicho:98b}) and we obtain a degree 9 equation. \[ {\textrm{Z}}: s^7 t^2 - 2 s^6 t^3 + s^5 t^4 + 9 s^5 t^2 u^2 + 9 s^4 t^2 v^3 - 6 s^4 t u^4 - 12 s^3 t u^2 v^3 + s^3 u^6 - 6 s^2 t v^6 + 3 s^2 u^4 v^3 + 3 s u^2 v^6 + v^9 \] We have that $p$ is birational on the hyperplane sections. We confirm that the sectional genus of ${\textrm{Z}}$ is endeed 4: Let ${\textrm{Y}}_e$ be the affine chart of ${\textrm{Y}}$  such that  $e\neq0$. Similarly for ${\textrm{Z}}_v$ and ${\textrm{Q}}_v$. Let $I({\textrm{Y}}_e)$ be the ideal of ${\textrm{Y}}_e$. Let ${\textbf{Q}}({\textrm{Z}}_v)$ be the function field of ${\textrm{Z}}_v$.

 We invert the projection $\ensuremath{{\textrm{Y}}\stackrel{p}{\rightarrow}{\textrm{Z}}}$ by considering the induced affine morphism on the charts $\ensuremath{{\textrm{Y}}_e\stackrel{p_{ev}}{\rightarrow}{\textrm{Z}}_v}$. From \cite{Schicho:98b} it follows that we can compute the inverse of $p_{ev}$ by computing the Gr\"obner basis of $I({\textrm{Y}}_e)\overset{}{\underset{{\textbf{Q}}}{\otimes}}{\textbf{Q}}({\textrm{Z}}_v) + {\textbf{Q}}({\textrm{Z}}_v)\langle s-a,t-d, g-u\rangle $. Composing $\ensuremath{{\textrm{Z}}\stackrel{p^{-1}}{\dashrightarrow}{\textrm{Y}}}$ with $\ensuremath{{\textrm{Y}}\stackrel{q}{\rightarrow}{\textrm{Q}}}$ we obtain a 2:1 morphism: \[ \ensuremath{r:{\textrm{Z}}\rightarrow {\textrm{Q}},\quad  \rho=(s:t:u:v) \mapsto  (r_0(\rho):r_1(\rho):r_2(\rho):r_3(\rho)) }. \] Let $S_5$ be the ideal in Example~\ref{ex:f2_dual_tangent_developable}.

 If $(k:l:m:n) \in V(S_5)$ then $ks+lt+mu+nv=0$ is a bitangent plane of ${\textrm{C}}$. From Proposition~\ref{prop:f2_dp1_D} it follows that the bitangent planes sections pull back along $r$ to a family of rational curves on ${\textrm{Z}}$.

 We parametrize $V(S_5)$: \[ \ensuremath{\gamma:{\textbf{P}}\rightarrow V(S_5),\quad \alpha=(\alpha_0:\alpha_1)\mapsto (\gamma_0(\alpha):\ldots:\gamma_3(\alpha))=(27\alpha_1^2\alpha_0:4\alpha_0^3:-27\alpha_1^3:27\alpha_1^3)}, \] and obtain a family of bitangent planes indexed by the projective line: $ \gamma_0 s+\gamma_1 t + \gamma_2 u + \gamma_3 v $.

 We compute the pull back of this family of bitangent plane sections by considering the Gr\"obner basis of ${\textbf{Q}}\langle \gamma_0 r_0+\gamma_1 r_1 + \gamma_2 r_2 + \gamma_3 r_3\rangle  + I({\textrm{Z}})$. It follows that a chart of the corresponding family of curves on ${\textrm{Z}}_v$ is defined by the moving components of ${\textrm{Z}}_v\cap {\textrm{H}}_\kappa$ with \[ {\textrm{H}}_\kappa: -12 s^3 t \kappa^3 + 4 s^2 u^2 \kappa^3 + 27 s^3 t - 27 s^2 t^2 - 81 s^2 t \kappa + 27 s u^2 \kappa + 4 s \kappa^3 + 27 \kappa, \] where $\kappa=\frac{\alpha_0}{\alpha_1} \in {\textbf{Q}}$ with $\alpha_1\neq 0$. We computed ${\textrm{H}}_\kappa$ by considering the irreducible factor of the shortest equation in this Gr\"obner basis, which depends on $\kappa$. Let $F\subset {\textrm{Z}}\times{\textbf{P}}$ be the family of curves defined by ${\textrm{Z}}_v\cap {\textrm{H}}_\kappa$.

 We compute the projection of ${\textrm{H}}_\kappa\cap{\textrm{Z}}_v$ from a point outside ${\textrm{Z}}_v$ to the $u,t$ plane by eliminating $s$ from the ideal $I({\textrm{Z}}_v)+I({\textrm{H}}_\kappa)$. We factor the resulting equation and consider the irreducible equation of a curve in the $u,t$ plane which depends on $\kappa$. We confirm that this defines a rational curve of degree 6.

 The divisor class of a curve in $F$ is $-2K$. We have that $D(-2K)=6$, so we indeed expect a curve of degree 6. From \cite{nls1} it follows that $-KF\geq 2$ and thus $F$ defines indeed a minimal family.

 We want to determine the number $\lambda$ of curves in $F$ which go to a generic point on ${\textrm{Z}}$. Similar as in the example of the non-fibration family defined by lines tangent to the unit circle we find that $\lambda$ equals the mapping degree of the unirational parametrization defined by $F$.

 We have that $\lambda$ equals the number of bitangent planes of ${\textrm{C}}$ in $V(S_5)$ through a generic point $\tau$ on the quadric cone ${\textrm{Q}}$ outside ${\textrm{C}}$.

 We project ${\textrm{C}}$ from this point $\tau$ and obtain a curve ${\textrm{N}} \subset{\textbf{P}}^2$ of degree 6 and geometric genus 1.

 The bitangent lines of ${\textrm{N}}$ pull back to bitangent planes of ${\textrm{C}}$ through $\tau$. Let ${\textrm{N}}^*\subset{\textbf{P}}^2$ be the dual curve of ${\textrm{N}}\subset{\textbf{P}}^2$. Thus a point on ${\textrm{N}}^*$ corresponds to a tangent line of ${\textrm{N}}$.

 From \cite{har1}, chapter 4, section 3, it follows that projecting space curves to the plane only introduces $A_1$ singularities. From the Pl\"ucker formulas for plane curves (see appendix B, section 6, proposition 174 in \cite{nls2} or \cite{jha2}) it follows that ${\textrm{N}}$ has $3A_2$ and $6A_1$ singularities. We find that ${\textrm{N}}^*$ has $15A_1$ and $12A_2$ singularities.

 The $15A_1$ singularities of ${\textrm{N}}^*$ correspond to bitangent lines of ${\textrm{N}}$. The Pl\"ucker formulas assume that the space curve is generic. Since ${\textrm{C}}$ lies on ${\textrm{Q}}$ we have tritangent planes, which do not occur for generic space curves. The tritangent plane through the vertex is projected to a tritangent line of ${\textrm{N}}$. The tritangent line of ${\textrm{N}}$ corresponds to the intersection of three branches on ${\textrm{N}}^*$. After some analytic deformation of ${\textrm{C}}$ this corresponds to $3A_1$  singularities on ${\textrm{N}}^*$.

 From Proposition~\ref{prop:f2_dp1_D} it follows that the tritangent plane through the vertex pulls back along $\ensuremath{{\textrm{Z}}\stackrel{r}{\rightarrow}{\textrm{Q}}}$ to an elliptic curve in $|D|$ and thus does not define not a curve in the family $F$. It follows $\lambda=15-3=12$.

 Since 12 different curves of $F$ go a generic point on ${\textrm{Z}}$, it follows that $F$ can not be defined by the fibres of a map, so $F$ is indeed a minimal non-fibration family.

\end{example}

\newpage
\section{Acknowledgements}
This paper is a part of the authors PhD thesis. It is my pleasure to acknowledge that the many discussions with my advisor Josef Schicho is a major contribution to this paper. The discovery of minimal non-fibration families should be attributed to him.

 I would like to thank Helmut Pottmann for reviewing this paper and suggesting some improvements.

 The images were made using Surfex (\cite{surfex}). The algorithms were implemented using the computer algebra system Sage  (\cite{sage}) with additional interfacing to Magma (\cite{magma}).

 This research was supported by the Austrian Science Fund (FWF): project P21461. \bibliography{geometry,schicho,rational_curves}\paragraph{Address of author:} ~\\
~\\
 King Abdullah University of Science and Technology, Thuwal, Kingdom of Saudi Arabia \\
 \textbf{email:} niels.lubbes@gmail.com\printindex
\end{document}